\documentclass[final]{siamltex}

\usepackage{epsfig}
\usepackage{amssymb}
\usepackage{amsmath}
\usepackage{amsfonts}
\usepackage[latin1]{inputenc}
\usepackage[english]{babel}

\title{Numerical study of a flow of regular planar curves that develop
singularities at finite time\thanks{This work was partly supported
by grant MTM2007-62186 of MEC (Spain). It was partly announced in
the proceedings of the XX CEDYA, at Sevilla, Spain (2007).}}

\author{Francisco de la Hoz\thanks{Departamento de Matemática
Aplicada // Escuela Universitaria de Ingeniería Técnica Industrial
// Universidad del País Vasco // Plaza de la Casilla, 3 // 48012 Bilbao
(Spain) ({\tt francisco.delahoz@ehu.es}).}}

\begin{document}

\maketitle

\begin{abstract}
In this paper, we will study the following geometric flow, obtained
by Goldstein and Petrich while considering the evolution of a vortex
patch in the plane under Euler's equations,
$$
\mathbf X_t = -k_s\mathbf n - \dfrac{1}{2}k^2\mathbf T,
$$

\noindent with $s$ being the arc-length parameter and $k$ the
curvature. Perelman and Vega proved in \cite{perelman} that this flow
has a one-parameter family of regular solutions that develop a
corner-shaped singularity at finite time. We will give a method to
reproduce numerically the evolution of those solutions, as well as the
formation of the corner, showing several properties associated to them.
\end{abstract}



\pagestyle{myheadings} \thispagestyle{plain} \markboth{Francisco de
la Hoz}{Numerical study of a flow that develops singularities at
finite time}

\section{Introduction}

In this paper, we will consider the following geometric flow of
planar curves that can develop singularities at finite time
\begin{equation}
\label{flujo} \mathbf X_t = -k_s\mathbf n - \dfrac{1}{2}k^2\mathbf
T,
\end{equation}

\noindent with $s$ being the arc-length parameter and $k$ the
curvature, $\mathbf T_s = k\mathbf n$. It was obtained by Goldstein and
Petrich in \cite{goldstein}; their motivation was the problem of the
evolution of a vortex patch in the plane subject to Euler equations
\cite{majda}. If the boundary is at least piecewise $\mathcal C^1$, the
exact motion of the boundary of the patch satisfies
\begin{align}
\label{parche} \mathbf X_t(\alpha, t) & =
-\frac{\omega_p}{2\pi}\int_0^{2\pi}\ln\left(\frac{\|\mathbf
X(\alpha, t) - \mathbf X(\alpha', t)\|}{r_0}\right)\mathbf
X_\alpha(\alpha', t)d\alpha',
\end{align}

\noindent where $\alpha$ is the Lagrangian parameter, $\omega_p$ is
the vorticity and $r_0$ is an arbitrary parameter whose choice does
not affect the dynamics of the curve.

Let us rewrite \eqref{parche} using the arc-length parameter $s$; then
\begin{align}
\mathbf X_t(s, t) = -\frac{\omega_p}{2\pi}\int_{s - L(t)/2}^{s +
L(t)/2}\ln\left(\frac{\|\mathbf X(s, t) - \mathbf X(s',
t)\|}{r_0}\right)\mathbf X_s(s', t)ds',
\end{align}

\noindent with $L(t)$ being the length at time $t$. We truncate this
last integral by introducing a cutoff at $s = \pm\Lambda$:
\begin{align}
\label{cutoff} \mathbf X_t(s, t) \approx -\frac{\omega_p}{2\pi}\int_{s -
\Lambda}^{s + \Lambda}\ln\left(\frac{\|\mathbf X(s, t) - \mathbf
X(s', t)\|}{r_0}\right)\mathbf X_s(s', t)ds'.
\end{align}

\noindent We expand $\mathbf X(s', t)$ and $\mathbf X_s(s', t)$ into
powers of $\Delta = s - s'$, bearing in mind $\mathbf X_s = \mathbf T$,
$\mathbf T_s = k\mathbf n$, $\mathbf n_s = -k\mathbf T$:
\begin{align*}
\mathbf X(s') & = \mathbf X(s) + \mathbf T(s)\Delta +
\frac{1}{2}k(s)\mathbf n(s)\Delta^2 + \cdots
\\
\mathbf T(s') & = \mathbf T(s) + k(s)\mathbf n(s)\Delta +
\frac{1}{2}[k_s(s)\mathbf n(s) - k^2(s)\mathbf
T(s)]\Delta^2 + \cdots
\end{align*}

\noindent Introducing the Taylor expansions into \eqref{cutoff}, we
can integrate term by term. If we represent the PDE for $\mathbf X$
as $\mathbf X_t = U\mathbf n + W\mathbf T$, then, considering the
leading terms of the expansions, we obtain the following
approximations for $U$ and $W$:
\begin{align*}
U & \approx
-\frac{\omega_p\Lambda^3}{6\pi}k_s(s)\left(\ln\Lambda -
\frac{1}{3} - \ln(r_0)\right),
\\
W & \approx -\frac{\omega_p\Lambda}{\pi}(\ln\Lambda - 1 - \ln(r_0))
+ \frac{\omega_p\Lambda^3}{6\pi}k^2(s)\left(\ln\Lambda -
\frac{1}{3} - \frac{1}{4} - \ln(r_0)\right).
\end{align*}

\noindent We choose $r_0 = \Lambda e^{-1/2}$ in order to have $W_s = k
U$ \cite{goldstein1}. With that choice, the final approximation for
$\mathbf X$ is
\begin{align}
\mathbf X_t \approx
-\frac{\omega_p\Lambda^3}{36\pi}\left[k_s\mathbf n +
\frac{1}{2}k^2\mathbf T\right] +
\frac{\omega_p\Lambda}{2\pi}\mathbf T.
\end{align}

\noindent A Galilean transformation removes the term
$\dfrac{\omega_p\Lambda}{2\pi}\mathbf T$; then, the factor
$-\dfrac{\omega_p\Lambda^3}{36\pi}$ is absorbed after a change of
variable, getting \eqref{flujo}
$$
\mathbf X_t = -k_s\mathbf n - \dfrac{1}{2}k^2\mathbf T.
$$

\noindent Since we are considering planar curves, we can identify
the plane where they live with $\mathbb C$; denoting $z\equiv\mathbf
X$ and bearing in mind that $z_{ss} = ikz_s$, the last equation
becomes
\begin{equation}
\label{z_t = -z_sss +...}
\begin{cases}
z_t = - z_{sss} + \dfrac{3}{2}\bar z_s z_{ss}^2, \\
|z_s|^2 = 1, \qquad t\not= 0.
\end{cases}
\end{equation}

\noindent In this form, the local induction approximation preserves
some of the basic conserved quantities of the exact vortex patch
dynamics; for example, area, center of mass and angular momentum
\cite{batchelor}. It also preserves $s_\alpha$ and, in particular, the
length, which is not true in the vortex patch problem. Indeed,
numerical calculations \cite{deem, dritschel} show that small bumps in
the boundary of isolated vortex patches of constant curvature cause
filamentation phenomena to occur, i.e., the ejection of thin filaments
into the surrounding fluid. Nevertheless, for a small enough
perturbation, the time at which filamentation appears can be made
arbitrarily large and we can assume that the initial parametrization of
the curve $s_\alpha$, and hence $L$, are time-independent.

Flow \eqref{z_t = -z_sss +...} is time-reversible, because if $z(s, t)$
is a solution, so is $z(-s, -t)$. It is completely determined by its
curvature, $k(s, t)$, except for a rigid movement that changes with
time and that can be fixed by the initial conditions. As shown by
Goldstein and Petrich, $k$ satisfies the modified Korteweg-De~Vries
(mKdV) equation
\begin{equation}
\label{mkdv} k_t + k_{sss} + {3\over2}k^2k_s = 0.
\end{equation}

\noindent To relate the vortex patch evolution and the mKdV equation,
Goldstein and Petrich followed previous ideas by Hasimoto
\cite{hasimoto}, who connected the nonlinear Schrödinger (NLS) equation
with the motion of vortex filaments in $\mathbb R^3$, ideas which were
extended by Lamb \cite{lamb}. Later, Nakayama, Segur and Wadati
\cite{wadati} identified the connection between integrable evolution
equations and the motion of curves in the plane and in $\mathbb R^3$.
More recently, Wexler and Dorsey \cite{wexler} found that under a local
induction approximation, the contour dynamics of the edge of a
two-dimensional electron system can be described again by the mKdV
equation.

In \cite{perelman}, Perelman and Vega proved the existence of a
regular family of solutions for \eqref{z_t = -z_sss +...} that
develop corner-shaped singularities at finite time. Conversely, they
also proved the existence of solutions of the mKdV equation
\eqref{mkdv} with initial conditions given by
\begin{equation}
\label{ks0} k(s, 0) = a\delta(s), \qquad a\in\mathbb R,
\end{equation}

\noindent where $\delta(s)$ is the Dirac delta function and $|a|$ is small enough. The
corresponding initial condition for \eqref{z_t = -z_sss +...} is
\begin{equation}
\label{zcorner} z(s, 0) =
\begin{cases}
z_0 + se^{i\theta^+}, & s \ge 0,
\\
z_0 + se^{i\theta^-}, & s \le 0,
\end{cases}
\end{equation}

\noindent for some $\theta^+$, $\theta^-\in[0, 2\pi)$, $\theta^+ - \theta^- = 2a$.

Perelman and Vega looked for self-similar solutions of \eqref{mkdv}
in the following form
\begin{equation}
\label{k autosemejante} k(s, t) =
{2\over(3t)^{1/3}}u\left({s\over(3t)^{1/3}}\right), \quad t > 0.
\end{equation}

\noindent which leads to study the following ODE
\begin{equation}
u_{xx} - xu + 2u^3 = \mu, \qquad x\in\mathbb R, \quad \mu\in\mathbb
R,
\end{equation}

\noindent being $\mu$ an integration constant. We will only consider
the case $\mu = 0$,
\begin{equation}
\label{ode u} u_{xx} - xu + 2u^3 = 0, \qquad x\in\mathbb R,
\end{equation}

\noindent but the method developed in this paper can be easily
implemented for $\mu\not=0$.

Equivalently, the self-similar solutions for \eqref{z_t = -z_sss
+...} are of the form
\begin{equation}
z(s, t) = t^{1/3}\omega\left({{s\over t^{1/3}}}\right), \qquad t >
0,
\end{equation}

\noindent which leads to study
\begin{equation}
\label{ode omega}
\begin{cases}
\dfrac{1}{3}\omega - \dfrac{s}{3}\omega_s = -\omega_{sss} +
\dfrac{3}{2}\bar\omega_s\omega_{ss}^2, \qquad s\in\mathbb R,
\\
|\omega_s|^2 = 1.
\end{cases}
\end{equation}

\noindent Bearing in mind all the previous arguments, Perelman and
Vega proved the following theorems:
\begin{theorem}
\label{teorema1} There is $\epsilon_0 > 0$ such that if $a^2 <
\epsilon_0$, then there exist $\theta^\pm\in [0, 2\pi)$ and $\omega$
an analytic solution of \eqref{ode omega} such that if
\begin{equation}
z(s, t) = t^{1/3}\omega\left({{s\over t^{1/3}}}\right), \qquad t >
0,
\end{equation}

\noindent then
\begin{enumerate}
\item[(i)] $z$ solves \eqref{z_t = -z_sss +...} for
$t > 0$ and
\begin{align*}
\left|z(s, t) - se^{i\theta^+}\chi_{[0,\infty)}(s) -
se^{i\theta^-}\chi_{(-\infty,0]}(s) \right| \le ct^{1 / 3};
\end{align*}

with $\chi$ being the characteristic function.

\item[(ii)] $\theta^+ - \theta^- = 2a$;

\item[(iii)] The curvature $k$ of $z$ satisfies
\eqref{k autosemejante} and \eqref{ode u}.

\end{enumerate}

\end{theorem}

\begin{theorem}
\label{teorema2} There is $\epsilon_0 > 0$ such that given any $a$,
with $a^2 < \epsilon_0$, there exists a bounded real analytic $u(x,
t)$ solution of \eqref{ode u} such that if
\begin{equation}
k(s, t) = {2\over(3t)^{1/3}}u\left({s\over(3t)^{1/3}}\right), \quad
t > 0,
\end{equation}

\noindent then
\begin{enumerate}
\item[(i)] $k$ solves the mKdV equation \eqref{mkdv};

\item[(ii)] $\displaystyle\int_{-\infty}^{+\infty}u(x)dx = a$;

\item[(iii)] $\displaystyle\lim_{t\to0}k(\cdot, t) = 2a\delta$ in
    $\mathcal S'$, the space of tempered distributions.

\end{enumerate}

\end{theorem}

\noindent These theorems, which constitute the theoretical basis of
this paper, guarantee the formation of a corner-shaped singularity at
finite time for \eqref{z_t = -z_sss +...}, provided that $|a|$ is small
enough. Nonetheless, numerical simulations in subsection
\ref{subsection experimentsth} will give evidence that singularity
formation happens also for any $a\in(-\pi/2, \pi/2)$, i.e., for
parameter values outside the scope of Perelman and Vega's theory.

The purpose of this work is to study the self-similar solutions of
\eqref{mkdv} from a numerical point of view, as well as the formation
of their corresponding corner \eqref{zcorner}, going backwards in time
from $t = 1$ until $t = 0$, because \eqref{mkdv} is also time
reversible.

Instead of developing a numerical method for \eqref{z_t = -z_sss +...}
or \eqref{mkdv}, we rather consider the angle
\begin{align}
\theta(s, t) = \theta(-\infty, t) + \int_{-\infty}^sk(s', t)ds',
\end{align}

\noindent and the PDE for the angle, obtained after integrating
\eqref{mkdv} once,
\begin{align}
\label{thetat} \theta_t(s, t) = -\theta_{sss}(s, t) -
\dfrac{1}{2}(\theta_s)^3(s, t).
\end{align}

\noindent Working with $\theta$ has two main advantages: it allows
to guarantee naturally $|z_s(s, t)| = 1$ for all $t$ and we can
preserve numerically the conserved quantity of \eqref{mkdv}
$$\dfrac{1}{2}\int_{-\infty}^{+\infty}u(s)ds = \int_{-\infty}^{+\infty}k(s,
t)ds = \theta^+ - \theta^-,$$

\noindent by fixing $\theta(-\infty, t) = \theta^-$ and
$\theta(+\infty, t) = \theta^+$, for all $t$.

The structure of this paper is as follows: In section \ref{sectionux},
we integrate \eqref{ode u}, looking for admissible initial data $(u(0),
u_x(0))$, such that the corresponding solutions satisfy $u(x) \to 0$ as
$x\to\infty$; this must be carefully done, because those solutions are
very unstable.

The admissible pairs $(u(0), u_x(0))$ form a connected curve. Each
point of this curve determines one solution for \eqref{z_t = -z_sss
+...}, \eqref{mkdv} and \eqref{thetat}; hence, we have one-parameter
families of regular solutions for \eqref{z_t = -z_sss +...},
\eqref{mkdv} and \eqref{thetat} that develop singularities at finite
time.

In section \ref{sectionmethod}, a spectral numerical method with
integrating factor for \eqref{thetat} is developed, having truncated
$s\in\mathbb R$ to $[s_a, s_b]$, with $s_a \ll -1$, $s_b\gg 1$. We
impose the boundary conditions $\theta(s_a, t) = \theta^-$,
$\theta(s_b, t) = \theta^+$, which is equivalent to fixing the tangent
vectors of $z$ at $s_a$ and $s_b$. In subsection \ref{subsection ths1}
we explain how to integrate $\theta(s, 1)$ in $[s_a, s_b]$, which
involves an estimate of $\int_{-\infty}^{+\infty}k(s, 1)$. Finally,
numerical experiments are carried out in subsection \ref{subsection
experimentsth}.

In the exact problem, the energy $\int_{-\infty}^{+\infty}k^2(s, t)ds =
\infty$ is a conserved quantity for all $t$; this infinite energy
concentrates at $s = 0$ as $t\to0$, which causes $z$ to develop the
corner-shaped singularity. In our numerical experiments, the energy in
$[s_a, s_b]$, $\int_{s_a}^{s_b}k^2(s, t)ds$, is finite, but,
nevertheless, it keeps approximately constant and it also tends to
concentrate at $s = 0$, as $t\to 0$. This fact shows that even after
having truncated $\mathbb R$ to $[s_a, s_b]$, the energy accumulation
process continues to be stable. It is also remarkable the good accuracy
with which we recover $k(0, t) = t^{-1/3}k(0, 1)$ even for small $t$,
hence approaching the Dirac delta function \eqref{ks0}. The numerical
results suggest that the accuracy of $k(0, t)$ could be improved
arbitrarily by increasing the length of $[s_a, s_b]$, i.e., the energy
of the system; it would be very interesting to prove analytically that
we can recover the solution of the exact problem by making
$s_a\to-\infty$ and $s_b\to\infty$.

In section \ref{sectionintegral}, we calculate the estimates of
$\int_{-\infty}^{+\infty}k(s, 1)$, as explained in subsection
\ref{subsection ths1}, for a large set of admissible initial data
$(u(0), u_x(0))$ of \eqref{ode u}, giving numerical evidence that
\begin{align*}
\int_{-\infty}^{+\infty} u(x)dx \in\left(-\dfrac{\pi}{2},
\dfrac{\pi}{2}\right) \mbox{ or, equivalently, }
\int_{-\infty}^{+\infty} k(s', t)ds'\in(-\pi, \pi), \quad\forall t.
\end{align*}

\noindent In section \ref{sectionclosed}, motivated by the original
vortex patch problem, we append to the initial datum $\theta(s, 1)$
a smooth function in such a way that the corresponding $z(s, 1)$ is
a closed regular curve without intersections. Numerical experiments
with this new initial datum are carried out in subsection
\ref{subsection experimentsthclosed}, showing that the method
developed in section \ref{sectionmethod} keeps $z(s, t)$ closed for
all $t$, preserving its inner area as well. It is also observed that
closing $z$ has no effect on the energy concentration process.

Finally, conclusions are summarized in section
\ref{sectionconclusions}.

Besides the local induction approximation done by Goldstein and Petrich
\cite{goldstein}, some other simplified models have been proposed to
describe the vortex patch dynamics. In \cite{constantin}, Constantin
and Titi introduced a hierarchy of area-preserving nonlinear
approximate equations, showing that the first of these equations,
starting from arbitrarily small neighborhoods of the circular vortex
patch, blows up. Later on, Alinhac \cite{alinhac} considered a
quadratic non area-preserving approximation for vortex patches with
contours near the unity circle, obtaining an instability result at
finite time.

In the vortex patch problem, Chemin has proven in \cite{chemin} that,
when considering a smooth initial contour, no finite-time singularities
may happen (infinite length, corners or cusps, for instance), i.e.,
smooth contours stay smooth in time; later on, this result has been
proven also by Bertozzi and Constantin \cite{bertozzi}. There is no
contradiction between those results and our model, because Kenig, Ponce
and Vega \cite{kenig} have proven in $\mathbb R$ that finite energy
solutions of the mKdV equation do not develop singularities; an
equivalent of this has been proved by Bourgain \cite{bourgain} in the
torus. In \cite{perelman}, Perelman and Vega were instead considering
infinite energy solutions of the mKdV equation, which caused the
singularity to happen.

Our experiments exhibit the energy concentration process in \eqref{z_t
= -z_sss +...}, both after having truncated $z$ and after having closed
it with a big loop. Since we consider finite energy solutions, we
cannot reproduce the corner-shaped singularity formation, but just
approach it. It would be very interesting to simulate the evolution of
this closed curve under the equations for the vortex patch. Although no
singularities can be obtained, it could be still expected to see some
energy concentration process taking place (obtaining ``big curvatures''
from ``small'' initial curvatures).

Despite its good properties (conservation of area, center of mass,
angular momentum) and the fact that it does not create
singularities, the main criticism of our model is that the
arc-length parameter is preserved for a simple closed curve and,
therefore, the total length of the curve remains constant. This does
not happen in the vortex patch problem, where there are examples in
which the length or curvature of the vortex patch boundary grow
rapidly \cite{alinhac, constantin}. Therefore, the mKdV equation
should be complemented by another equation for the evolution of
$s_\alpha$.

\section{Integration of $u(x)$}

\label{sectionux}

Let us come back to \eqref{ode u}
\begin{equation}
\label{airymodmu0} u_{xx} = xu - 2u^3, \qquad x\in\mathbb R.
\end{equation}

\noindent In order to integrate this second-order ODE, we rewrite it
as
\begin{equation}
\label{uxvx}
\begin{cases}
u_x = v, \\
v_x = xu - 2u^3
\end{cases}
\end{equation}

\noindent and use the fourth-order Runge-Kutta method, needing two
initial data $u(0)$ and $u_x(0) = v(0)$. $\gamma(x) = \int_0^xu(x')dx'$
can be also integrate, without extra effort, by adding the equation
$\gamma_x = u$ to the system \eqref{uxvx}. This idea will be used in
subsection \ref{subsection ths1}.

If $u$ is a solution of \eqref{airymodmu0}, then, for $x \ll 0$,
$u(x)$ has the same oscillatory behavior as the solutions of the
Airy equation $ u_{xx} - xu = 0$. Likewise, when $x\to\infty$, the
solutions of \eqref{airymodmu0} are characterized by being very
sensitive with respect to small variations of the initial data, i.
e., if $u(x)\to\infty$, small changes in the initial data can make
$u(x)\to-\infty$ and vice versa. Let us consider, for instance,
$u(0) = 0.024$, as well as four possible choices for $u_x(0)$ and
let us integrate numerically \eqref{airymodmu0}.
\begin{figure}[!ht]
\centering \resizebox{9cm}{4cm}{\includegraphics{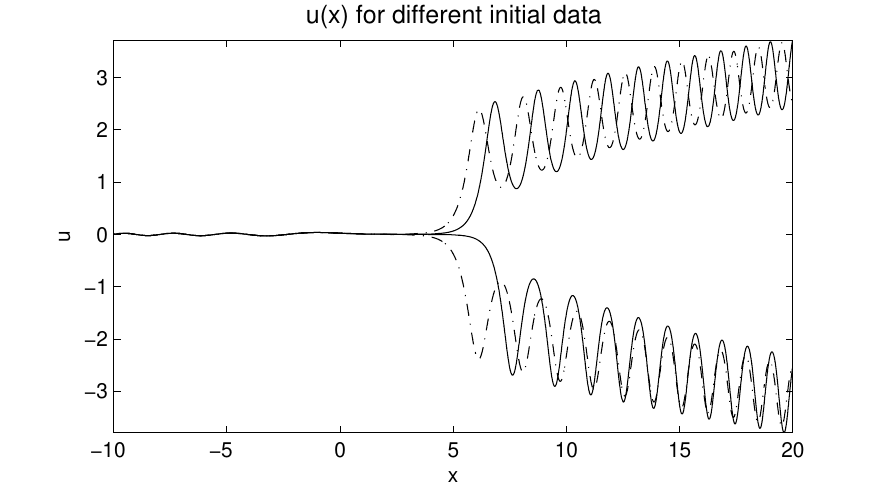}}
\caption{Dependance of $u(x)$ on the initial
conditions}\label{figura u(x) sensitiva}
\end{figure}

In figure \ref{figura u(x) sensitiva}, the solutions for $u_x(0) =
-0.018$ and $u_x(0) = -0.017$, plotted discontinuously, tend
respectively to $-\infty$ and $\infty$, and the solutions for
$u_x(0) = -0.0175$ and $u_x(0) = -0.0174$, plotted continuously,
tend respectively to $-\infty$ and $\infty$, although we observe
that for these last two values, the explosion happens a bit later.

If we go farther with the process, between $u_x(0) = -0.0175$ and
$u_x(0) = -0.0174$ there exists a unique $u_x(0)$ such that
$\lim_{x\to\infty}u(x) = 0$. That value, obtained with a bisection
technique, is approximately $u_x(0) = -0.0174881944\ldots$ In figure
\ref{figura u(x) sensitiva}, for scale reasons, we cannot distinguish
clearly the oscillations in the real negative axis; in figure
\ref{figura u(x) limite}, having used the limit value for $u_x(0)$,
those oscillations are displayed.
\begin{figure}[!ht]
\centering \resizebox{9cm}{4cm}{\includegraphics{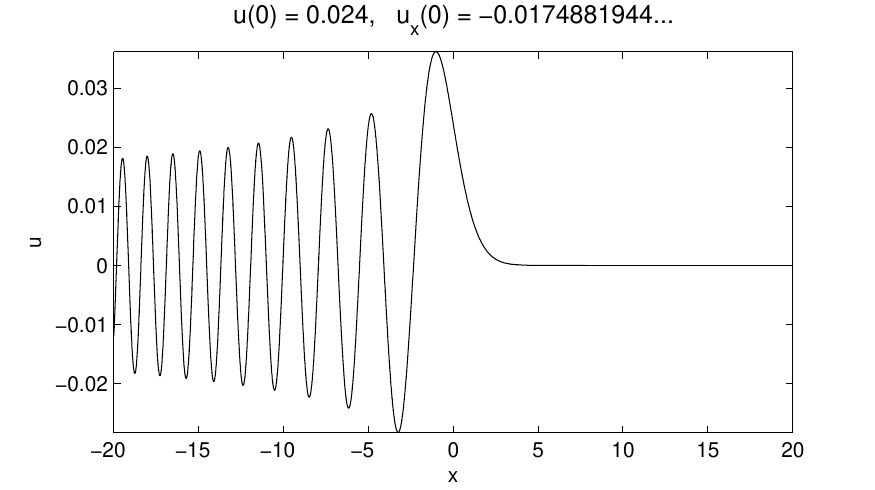}}
\caption{$u(x)$, with $\displaystyle{\lim_{x\to\infty}} u(x) = 0$.} \label{figura
u(x) limite}
\end{figure}

Because of the sensitivity with respect to initial data, we can
calculate numerically only as many decimals of the correct $u_x(0)$
as the machine precision allows us, which implies that we are only
delaying the explosion time for $u(x)$. We can say in an equivalent
way that the solutions that make
\begin{equation}
\label{limu(x)=0}\lim_{x\to\infty}u(x) = 0
\end{equation}

\noindent are highly unstable. Nonetheless, for our purposes, it is
enough to consider $u(x) \equiv 0$ for a big enough $x$, because
$u(x)\to0$ exponentially as $x\to\infty$. Hence, since
$|\int_{-\infty}^0u(x)dx| < \infty$, due to the fact that the
oscillations in $x < 0$ cancel one another, condition \eqref{limu(x)=0}
is equivalent to $|\int_{-\infty}^{+\infty}u(x)dx| < \infty$.
\begin{figure}[!ht]
\centering \resizebox{9cm}{6cm}{\includegraphics{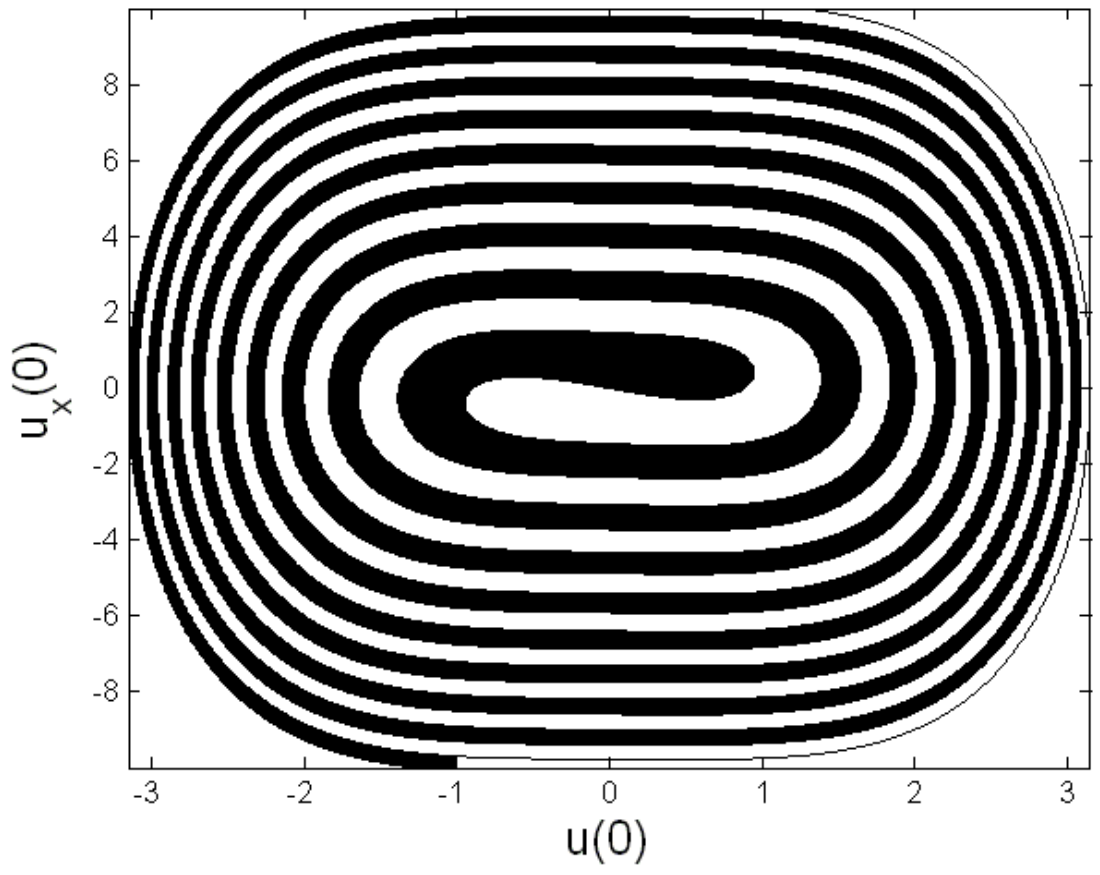}}
\caption{Pairs $(u(0), u_x(0))$, with
$\displaystyle{\lim_{x\to\infty}} u(x)=+\infty$ (black) and
$\displaystyle{\lim_{x\to\infty}} u(x)=-\infty$ (white)}
\label{figura u(x) dos regiones}
\end{figure}

In figure \ref{figura u(x) dos regiones}, the points belonging to the
region in black, when taken as initial data $(u(0), u_x(0))$ of
\eqref{airymodmu0}, give us solutions for \eqref{airymodmu0} such that
$u(x)\to+\infty$, as $x\to\infty$; for the points of the region in
white, we have $u(x)\to-\infty$ instead. The boundary (figure
\ref{figura u(x) pares admisibles}) between both regions are the
admissible pairs $(u(0), u_x(0))$ whose corresponding solutions satisfy
\eqref{limu(x)=0}. That boundary is a connected unbounded
one-dimensional curve that divides the plane in two antisymmetrical
halves. Each of the points of that curve determines one single
self-similar $k(s, t)$, one single self-similar $\theta(s, t)$ and one
single self-similar $z(s, t)$, which are solutions of \eqref{mkdv},
\eqref{thetat} and \eqref{z_t = -z_sss +...}, respectively. Therefore,
we have one-parameter families of solutions for \eqref{mkdv},
\eqref{thetat} and \eqref{z_t = -z_sss +...} that develop a singularity
at finite time.
\begin{figure}[!ht]
\centering \resizebox{9cm}{6cm}{\includegraphics{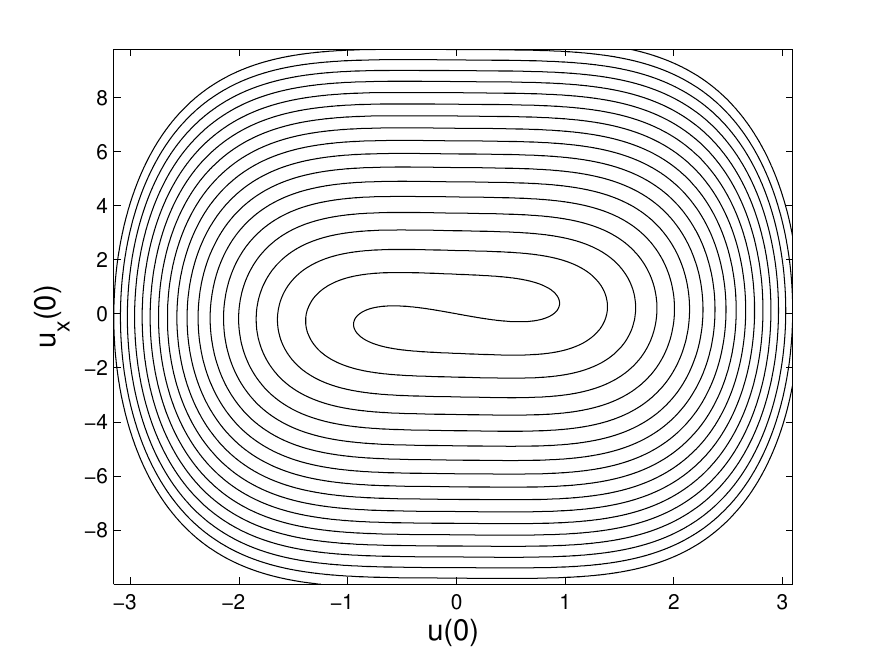}}
\caption{Admissible pairs. The points of this curve, when taken as
initial data of \eqref{uxvx}, make
$\displaystyle{\lim_{x\to\infty}}u(x) = 0$.} \label{figura u(x)
pares admisibles}
\end{figure}

\section{A numerical method for the evolution of \eqref{z_t = -z_sss
+...}}

\label{sectionmethod}

The $k(s, t)$, $\theta(s, t)$ and $z(s, t)$ we are considering are
characterized by
\begin{align*}
k(-\infty, t) & = 0, & k(+\infty, t) & = 0, \\
\theta(-\infty, t) & = \theta^-, & \theta(+\infty, t) & = \theta^+, \\
\lim_{s\to-\infty}|z(s, t) & - se^{i\theta^-}| = 0, &
\lim_{s\to+\infty}|z(s, t) & - se^{i\theta^-}| = 0.
\end{align*}

\noindent Taking $\theta(s, t)$ as the evolution variable, we can
preserver naturally $\int_{-\infty}^{+\infty}k(s, t)ds$ and $|z_s(s,
t)| = 1$. The PDE for $\theta$ with boundary conditions is
\begin{equation}
\begin{cases}
\theta_t(s, t) = -\theta_{sss}(s, t) - \dfrac{1}{2}(\theta_s)^3(s,
t), \qquad s\in\mathbb R, \\
\theta(-\infty, t) = \theta^-, \\
\theta(+\infty, t) = \theta^+.
\end{cases}
\end{equation}

\noindent Due to the difficulty of considering the whole $\mathbb
R$, we will restrict ourselves to $s\in[s_a, s_b]$, with $s_a \ll
-1$ and $s_b \gg 1$
\begin{equation}
\label{theta_t[s_a,s_b]}
\begin{cases}
\theta_t(s, t) = -\theta_{sss}(s, t) - \dfrac{1}{2}(\theta_s)^3(s,
t), \qquad s\in[s_a, s_b], \\
\theta(s_a, t) = \theta^-, \\
\theta(s_b, t) = \theta^+.
\end{cases}
\end{equation}

\noindent $\theta(s, t)$ is not periodic in $[s_a, s_b]$, so we
define
\begin{equation}
\label{theta periodica t}\tilde\theta(s, t) = \theta(s, t) - {s -
s_a\over L}(\theta^+ - \theta^-) - \theta^-, \qquad L = s_b - s_a,
\end{equation}

\noindent which is periodic and regular for all $t > 0$.
(\ref{theta_t[s_a,s_b]}) gets transformed into
\begin{equation}
\label{tildetheta_t[s_a,s_b]}
\begin{cases}
\tilde\theta_t(s, t) = -\tilde\theta_{sss}(s, t) -
\dfrac{1}{2}\left(\tilde\theta_s(s, t) +
\dfrac{\theta^+-\theta^-}{L}\right)^3, \qquad s\in[s_a, s_b], \\
\tilde\theta(s_a, t) = \tilde\theta(s_b, t) = 0.
\end{cases}
\end{equation}

\noindent We will calculate the numerical evolution of
$\tilde\theta$ at $N + 1$ equidistant points $\{s_j\}$ in $[s_a,
s_b]$. Denoting with some abuse of notation $\hat\theta(\xi, t)
\equiv (\tilde\theta(\xi, t))^\wedge$, we represent
\begin{equation}
\tilde\theta(s_j, t) = \sum_{\xi = -N/2}^{N/2-1}\hat \theta(\xi,
t)\exp\left[{2\pi i\xi\over L}(s_j - s_a)\right],
\end{equation}

\noindent where $L = s_b - s_a$, $\Delta s = \frac{L}{N}$ and $s_j =
s_a + j\Delta s$. Therefore, (\ref{tildetheta_t[s_a,s_b]}) is
transformed into
\begin{equation}
\hat\theta_t(\xi, t) = -\left(\dfrac{2\pi
i\xi}{L}\right)^3\hat\theta(\xi, t) -
\left[\dfrac{1}{2}\left(\tilde\theta_s(s, t) +
\dfrac{\theta^+-\theta^-}{L}\right)^3\right]^{\bigwedge} (\xi, t),
\end{equation}

\noindent with $\xi = -\frac{N}{2}, \cdots, \frac{N}{2} - 1$. Since
we are working with the frequency, the third derivative is
transformed into a multiplier that can be absorbed by means of an
integrating factor
$$
\left\{\hat\theta(\xi, t)\exp\left[t\left(\dfrac{2\pi
i\xi}{L}\right)^3\right]\right\}_t \!\!\! =
-\exp\left[t\left(\dfrac{2\pi i\xi}{L}\right)^3\right]\!\!\!
\left[\dfrac{1}{2}\left(\tilde\theta_s(s, t) +
\dfrac{\theta^+-\theta^-}{L}\right)^3\right]^{\bigwedge} \!\!\!
(\xi, t).
$$

\noindent The advantage of using an integrating factor is twofold:
it allows to integrate exactly the linear part of
\eqref{theta_t[s_a,s_b]}, increasing the accuracy of the numerical
results, and it relaxes considerably the time-step restrictions.

We apply the fourth-order Runge-Kutta in time with integrating
factor as described in \cite{spectralmatlab}. Denoting
\begin{align*}
\tilde\Theta^0(s) & = \tilde\theta(s, t^0), & \hat\Theta^0(\xi) & =
\hat\theta(\xi, t^0),
\\
\tilde\Theta^n(s) & \approx\tilde\theta^n(s) \equiv \tilde \theta(s,
t^n), & \hat\Theta^n(\xi) & \approx\hat\theta^n(\xi) \equiv
\hat\theta(\xi, t^n), \qquad t^n = t^0 + n\Delta t,
\end{align*}

\noindent we have
\begin{align*}
\hat A(\xi) & = -\dfrac{1}{2}\left\{\left(\left[\left(\dfrac{2\pi
i\xi}{L}\right)\hat\Theta^n(\xi) \right]^{\bigvee}+
\dfrac{\theta^+-\theta^-}{L}\right)^3\right\}^{\bigwedge}(\xi),
\\
\hat\Theta^{(A)}(\xi) & = \exp\left[-\dfrac{\Delta
t}{2}\left(\dfrac{2\pi
i\xi}{L}\right)^3\right]\left(\hat\Theta^n(\xi) + \dfrac{\Delta
t}{2}\hat A(\xi)\right),
\\
\hat B(\xi) & = -\dfrac{1}{2}\left\{\left(\left[\left(\dfrac{2\pi
i\xi}{L}\right)\hat\Theta^{(A)}(\xi) \right]^{\bigvee}+
\dfrac{\theta^+-\theta^-}{L}\right)^3\right\}^{\bigwedge}(\xi),
\\
\hat\Theta^{(B)}(\xi) & = \exp\left[-\dfrac{\Delta
t}{2}\left(\dfrac{2\pi i\xi}{L}\right)^3\right] \hat\Theta^n(\xi) +
\dfrac{\Delta t}{2}\hat B(\xi),
\\
\hat C(\xi) & = -\dfrac{1}{2}\left\{\left(\left[\left(\dfrac{2\pi
i\xi}{L}\right)\hat\Theta^{(B)}(\xi) \right]^{\bigvee}+
\dfrac{\theta^+-\theta^-}{L}\right)^3\right\}^{\bigwedge}(\xi),
\end{align*}

\begin{align*}
\hat\Theta^{(C)}(\xi) & = \exp\left[-\Delta t\left(\dfrac{2\pi
i\xi}{L}\right)^3\right]\hat\Theta^n(\xi) + \Delta
t\exp\left[-\dfrac{\Delta t}{2}\left(\dfrac{2\pi
i\xi}{L}\right)^3\right]\hat C(\xi),
\\
\hat D(\xi) & = -\dfrac{1}{2}\left\{\left(\left[\left(\dfrac{2\pi
i\xi}{L}\right)\hat\Theta^{(C)}(\xi) \right]^{\bigvee}+
\dfrac{\theta^+-\theta^-}{L}\right)^3\right\}^{\bigwedge}(\xi),
\\
\hat\Theta^{n + 1}(\xi) & = \exp\left[-\Delta t\left(\dfrac{2\pi
i\xi}{L}\right)^3\right]\hat\Theta^n(\xi) + \dfrac{\Delta
t}{6}\Bigg\{\exp\left[-\Delta t\left(\dfrac{2\pi
i\xi}{L}\right)^3\right]\hat A(\xi) \\
& + 2\exp\left[-\dfrac{\Delta t}{2}\left(\dfrac{2\pi
i\xi}{L}\right)^3\right]\left[\hat B(\xi) + \hat C(\xi)\right] +
\hat D(\xi)\Bigg\},
\\
t^{n + 1} & = t^n + \Delta t.
\end{align*}

\noindent The symbols $\wedge$ and $\vee$ denote respectively the
direct and inverse fast Fourier transforms (FFT) \cite{FFTW05}.

Finally, we force in every time step
\begin{align}
\begin{cases}
\tilde \Theta^{n + 1} = (\hat\Theta^{n + 1})^\vee
\equiv\Re((\hat\Theta^{n + 1})^\vee)
\\
\tilde \Theta^{n + 1}(s_0) = \tilde \Theta^{n + 1}(s_N) \equiv 0.
\end{cases}
\end{align}

\noindent Rounding $\tilde \Theta^{n + 1}$ to zero at the boundary
points $s_a = s_0$ and $s_b = s_N$ avoids the accumulation of little
errors of order $\mathcal O(10^{-9})$.

\subsection{Computation of $\theta(s, 1)$}

\label{subsection ths1}

Bearing in mind \eqref{k autosemejante}, at $t = 1$, $\theta(s, 1)$
is given by
\begin{align}
\label{thetas1} \theta(s, 1) - \theta(0, 1) = \int_0^sk(s', 1)ds' =
{2\over3^{1/3}}\int_0^s u\left({s'\over3^{1/3}}\right)ds' =
2\gamma\left({s\over3^{1/3}}\right),
\end{align}

\noindent where
\begin{align}
\gamma(x) = \int_0^xu(x')dx';
\end{align}

\noindent hence, $\gamma_x = u$ and \eqref{uxvx} can be generalized
to
\begin{equation}
\label{uxvxgx}
\begin{cases}
\gamma_x = u, \\
u_x = v, \\
v_x = xu - 2u^3,
\end{cases}
\end{equation}

\noindent with initial data $u(0)$, $u_x(0) = v(0)$ and $\gamma(0) =
0$. $(u(0), v(0))$ must be an admissible initial pair for \eqref{uxvx},
i.e., such that the corresponding $u(x)$ satisfies
$\lim_{x\to\infty}u(x) = 0$.

From now on, we will choose $u(0) = 0.72$, $v(0) = 1.1601860809647328$.
In figure \ref{figura theta(s,1)}, we have integrated \eqref{uxvxgx} in
$x\in[-80, 20]$, with $|\Delta x| = 10^{-5}$. For $x \ge 9.204$, we
have taken $k(x) \equiv 0$, because $|u(x)| < 10^{-12}$. We obtain
$\theta(s, 1)$ from \eqref{thetas1}; hence, $s = 3^{1/3}x$ and
$s\in[\tilde s_a, \tilde s_b] = [-115.38, 28.84]$. Since we have one
degree of freedom, we fix $\theta^+ = \theta(\tilde s_b, 1) \equiv 0$.

As we observe in the lower part of figure \ref{figura theta(s,1)},
$\theta(s, 1)$ tends to $\theta^-$ extremely slowly as $s\to-\infty$.
Nevertheless, we can determine $\lim_{s\to-\infty}\theta(s, 1)$ with
high accuracy, as the mean value of the first maximum and the first
minimum of $\theta(s, 1)$, for $s > \tilde s_a$. In our example, the
first maximum takes place at $s = -115.31$, being $\theta(-115.31, 1) =
-2.8403$; the first minimum takes place at $s = -114.80$, being
$\theta(-114.79, 1) = -3.0191$. Thus,
\begin{align*}
\lim_{s\to-\infty}\theta(s, 1) = \theta^-\approx-2.9297
\Rightarrow\int_{-\infty}^{+\infty}k(s', 1)ds' = \theta^+ - \theta^-
\approx 2.9294,
\end{align*}

\noindent quantity plotted in the lower part of figure \ref{figura
theta(s,1)} with a thinner stroke. In section \ref{sectionintegral},
we will approximate the value of $\int_{-\infty}^{+\infty}k(s',
1)ds'$ as a function of the admissible pairs $(u(0), u'(0))$.
\begin{figure}[!ht]
\centering \resizebox{9cm}{6cm}{\includegraphics{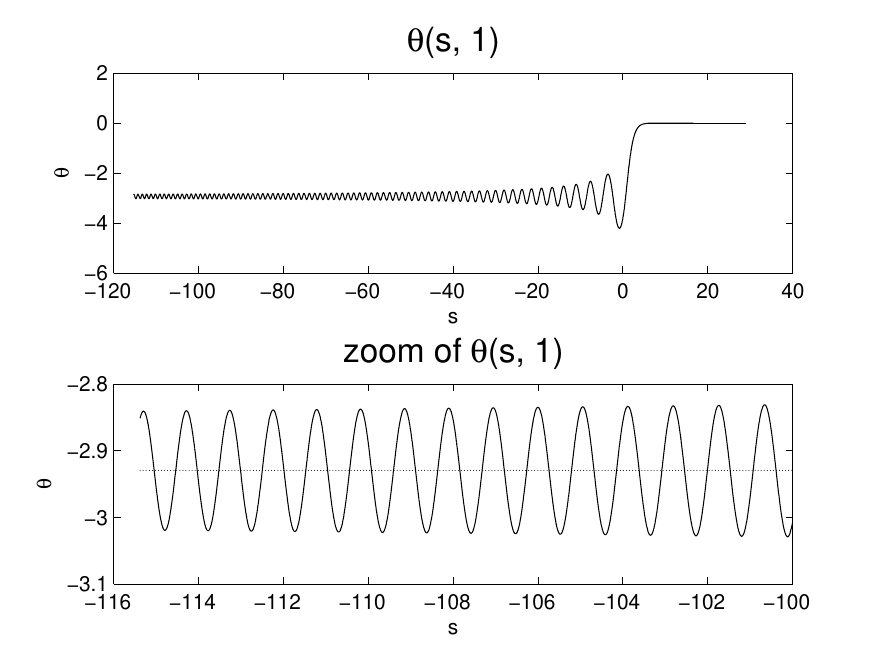}}
\caption{$\theta(s, 1)$. In the zoomed image, we have plotted with a
thinner stroke the mean between the first maximum and the first
minimum of $\theta(s, 1)$, for $s > s_a$, which gives a good estimate of $\theta^-$.}
\label{figura theta(s,1)}
\end{figure}

The $\theta(s, 1)$ obtained after integrating \eqref{uxvxgx} is not
directly suitable as an initial datum of \eqref{theta_t[s_a,s_b]},
because $\theta(\tilde s_a, 1)\not=\theta^-$ and (eq.~\ref{theta
periodica t})
\begin{equation}
\label{tildethetas1} \tilde\theta(s, 1) = \theta(s, 1) - {s -
\tilde s_a\over L}(\theta^+ - \theta^-) - \theta^-, \qquad L = \tilde s_b - \tilde s_a,
\end{equation}

\noindent is not periodic in $[\tilde s_a, \tilde s_b]$. Hence, we
have to modify $\theta(s, 1)$ as follows:
\begin{itemize}
\item We chose $N_1\in\mathbb N$ and divide $[\tilde s_a, \tilde
s_b]$ into $N_1 + 1$ equidistant points $\{s_j\}$, in such a way
that $0\in\{s_j\}$:
$$
s_j = \tilde s_a + j\Delta s, \quad \Delta s = \dfrac{\tilde s_b -
\tilde s_a}{N_1}, \qquad j = 0, \ldots, N_1.
$$
We have taken $N_1 = 2000$; hence $\Delta s\approx 0.07211$ and
$s_{1601} = 0$.

\item We find the $s_j$ with the lowest index, such that $s_j >
    s_{min}$, where $s_{min}$ is the first minimum point of
    $\theta(s, 1)$ in $(s_a, \infty)$, and $s_j$ satisfies that
    $\theta(s_j, 1)
    > \theta^-$ and $k(s_j, 1) = \theta_s(s_j, 1) > 0$. We name it
    $s_{joint}$. Here, $s_{joint} = s_{13}\approx -114.51$,
    $\theta(s_{joint}, 1) \approx -2.9150$ and $k(s_{joint}, 1)
    \approx 0.5479$.

\item We redefine $\theta(s, 1)$, appending at $s = s_{joint}$ a
smooth function that has a first order contact with $\theta(s, 1)$
and tends exponentially to $\theta^-$ for $s < \tilde s_a$. The
final expression for $\theta(s, 1)$ is

\begin{align*}
\theta(s, 1) =
\begin{cases}
\theta(s, 1), & s > s_{joint},
\\
\theta^- + [\theta(s_{joint}, 1) -
\theta^-]\exp\left[\dfrac{k(s_{joint}, 1)(s -
s_{joint})}{[\theta(s_{joint}, 1) - \theta^-]}\right], & s \le
s_{joint},
\end{cases}
\end{align*}

where, with some notational abuse, $\theta(s, 1)$ stands for the
original and the corrected functions.

\item We evaluate the new $\theta(s, 1)$ at the former $N_1 + 1$
points $\{s_j\}$, as well as in some new equidistant points
$\{s_j\}$ outside $[\tilde s_a, \tilde s_b]$. It is important that
the final number of $\{s_j\}$ is $N + 1 \equiv 2^n + 1$, for some
$n\in\mathbb N$, in order to apply FFT efficiently. In the new set
$\{s_j\}_0^N$, we define $s_a = s_0$, $s_b = s_N$.

\end{itemize}

\noindent The new function $\theta(s, 1)$ satisfies $\theta(s_a) =
\theta^-$, $\theta(s_b) = \theta^+$ and makes \eqref{theta periodica t}
periodic. We improve its regularity by applying a smooth spectral
filter to \eqref{theta periodica t}
\begin{equation}
\label{filtro}\hat \theta(\xi, 1) = \hat \theta(\xi,
1)\exp\left[-10\left({2.5|\xi| \over N}\right) ^ {25}\right],
\end{equation}

\noindent where $\hat\theta(\xi, 1) \equiv (\tilde\theta(\xi,
1))^\wedge$.

\begin{figure}[!ht]
\centering \resizebox{9cm}{6cm}{\includegraphics{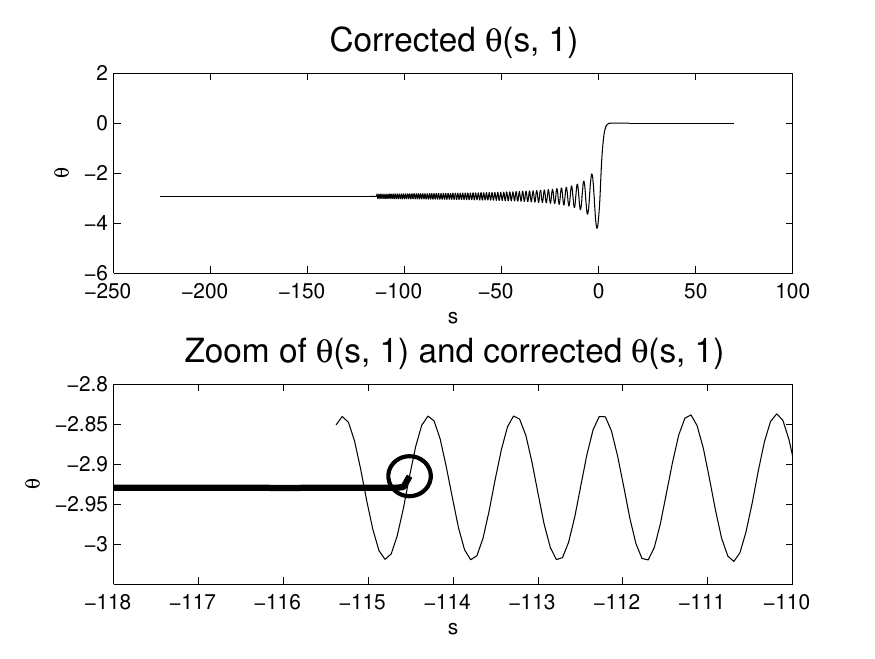}}
\caption{In the upper part, the final corrected $\theta(s, 1)$ is
plotted. In the lower part, both the original and the corrected
$\theta(s, 1)$ are plotted and zoomed near the joint (marked with a
circle); the corrected $\theta(s, 1)$ appears with a thicker
stroke.} \label{figura th corregida}
\end{figure}
In figure \ref{figura th corregida}, we have taken $N = 4096$. In the
upper part, we have plotted the final corrected $\theta(s, 1)$, with
$s\in[-226.14, 69.23]$. We observe two long constant segments at the
extremes of $\theta(s, 1)$, because of the rather large choice of $N$;
this is convenient to avoid periodicity phenomena, since the exact
$\theta(s, 1)$ is not periodic. In the lower zoomed part, we have
plotted both the original and the corrected $\theta(s, 1)$,
highlighting the point $s = s_{joint}$ with a circle.

In its new form, it is immediate to obtain $\theta_s(s, t)$ from
$\theta(s, t)$, through spectral derivation. Using again
representation \eqref{theta periodica t}, we get
\begin{equation}
\label{thetast} \theta_s(s, t) = \left[\theta(s, t) - {s - s_a\over
L}(\theta^+ - \theta^-) - \theta^-\right]_s + \dfrac{\theta^+ -
\theta^-}{L}, \qquad L = s_b - s_a.
\end{equation}

\noindent Although the process to obtain $\theta(s, 1)$ may look rather
artificial, if we calculate $\theta_s(s, 1)$ through \eqref{thetast}
and compare it with the $k(s, 1)$ obtained from integration of
\eqref{uxvxgx} and \eqref{k autosemejante}, we get the following error
table:
\begin{align*}
& \max_{s\in[-113.3, \infty)}|\theta_s(s, 1) - k(s, 1)| <
10^{-3}, \\
& \max_{s\in[-110.6, \infty)}|\theta_s(s, 1) - k(s, 1)| <
10^{-6}, \\
& \max_{s\in[-108.6, \infty)}|\theta_s(s, 1) - k(s, 1)| <
10^{-8}.
\end{align*}

\begin{figure}[!ht]
\centering \resizebox{9cm}{5cm}{\includegraphics{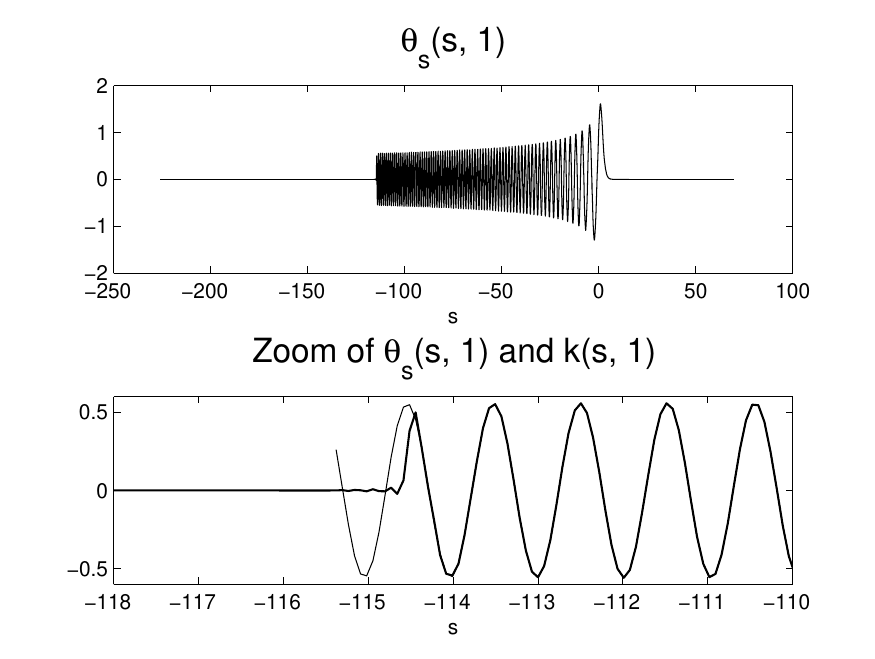}}%
\caption{In the upper part, the derivative of the final corrected
$\theta(s, 1)$ is plotted; observe that its support is finite. In
the lower part, the derivative of the corrected $\theta(s, 1)$,
plotted with a thicker stroke, is compared with the original $k(s,
1)$.} \label{figura ths vs k}
\end{figure}

\noindent Thus, $k(s, 1)$ is recovered from the corrected $\theta(s,
1)$ with high accuracy, except near the joint, $s_{joint} = -114.51$.
This is graphically illustrated in figure \ref{figura ths vs k};
observe that the support of $\theta_s(s, 1)$ is now finite.

$z(s, t)$ can also be immediately recovered from $\theta(s, t)$,
except for a rigid movement, because $z_s(s, t) = \exp(i\theta(s,
t))$. To fully determine $z(s, t)$, it is not complicate to see that
\begin{equation}
\label{z(0,t)0} z(0, t) = -2(3t)^{1/3} \left[iu'(0) +
u^2(0)\right]z_s(0, t),
\end{equation}

\noindent where $z_s(0, t)$ is a time-independent constant with unit
modulus.

\subsection{Numerical experiments}

\label{subsection experimentsth}

The aim of the method we have developed is to try to describe
numerically the formation of the singularity in the self-similar
solutions of equation \eqref{z_t = -z_sss +...} satisfying \eqref{k
autosemejante}. This corner-shaped singularity happens at finite, time;
indeed, at $t = 0$, we have
\begin{equation}
z(s, 0) =
\begin{cases}
z_0 + se^{i\theta^+}, & s\ge0, \\
z_0 + se^{i\theta^-}, & s\le0.
\end{cases}
\end{equation}

\noindent Equivalently, when $t \to 0$, the curvature tends to a
Dirac delta function
\begin{equation}
k(s, 0) = a\delta(s), \qquad a = \int_{-\infty}^{+\infty}u(x)dx.
\end{equation}

\noindent This happens because these self-similar solutions have
infinite energy
$$
\int_{-\infty}^{+\infty}k^2(s, t)ds = \infty, \qquad \forall t,
$$

\noindent and it tends to concentrate at $s = 0$ as $t \to 0$.

In our numerical experiments, we are not considering the whole
$\mathbb R$, but $s\in[s_a, s_b]$. At $s = s_a$ and $s = s_b$, we
fix the tangent vector of $z$, i.e., $\theta(s_a, t) = \theta^-$,
$\theta(s_b, t) = \theta^+$, for all $t$. As we observe in the upper
part of figure \ref{figura ths vs k}, the support of the initial
$k(s, 1)$ of our numerical experiments is finite; hence, our
numerical solutions have finite energy at $t = 1$:
$$
\int_{-\infty}^{+\infty}k^2(s, 1)ds = \int_{s_a}^{s_b}k^2(s, 1)ds <
\infty.
$$

\noindent We have executed the method with the initial $\theta(s,
1)$ plotted in figure \ref{figura th corregida}, i.e.,
\begin{align}
\label{experimento1}
\begin{cases}
u(0) = 0.72, \qquad u'(0) = 1.1601860809647328, \\
u(x) \mbox{ integrated when $x$} \in [-80, 9.203], \\
s \in [-226.14, 69.23], \qquad N = 4096, \qquad \Delta s = 0.07211.
\end{cases}
\end{align}

\noindent In order to measure the quality of our results, we analyze
two quantities: The evolution of the energy at $s\in[s_a,s_b]$
\begin{equation}
\label{Esasb} \int_{s_a}^{s_b}k^2(s, t)ds, \qquad t \in [0, 1],
\end{equation}

\noindent and the curvature at the origin $s = 0$, $k(0, t)$. From
\eqref{k autosemejante},
\begin{equation}
\label{k0t} k(0, t) = {2\over(3t)^{1/3}}u(0).
\end{equation}

\noindent When considering the non-truncated problem with
$s\in\mathbb R$, since all the infinite energy tends to concentrate
\begin{figure}[!ht]
\centering \resizebox{9cm}{5cm}{\includegraphics{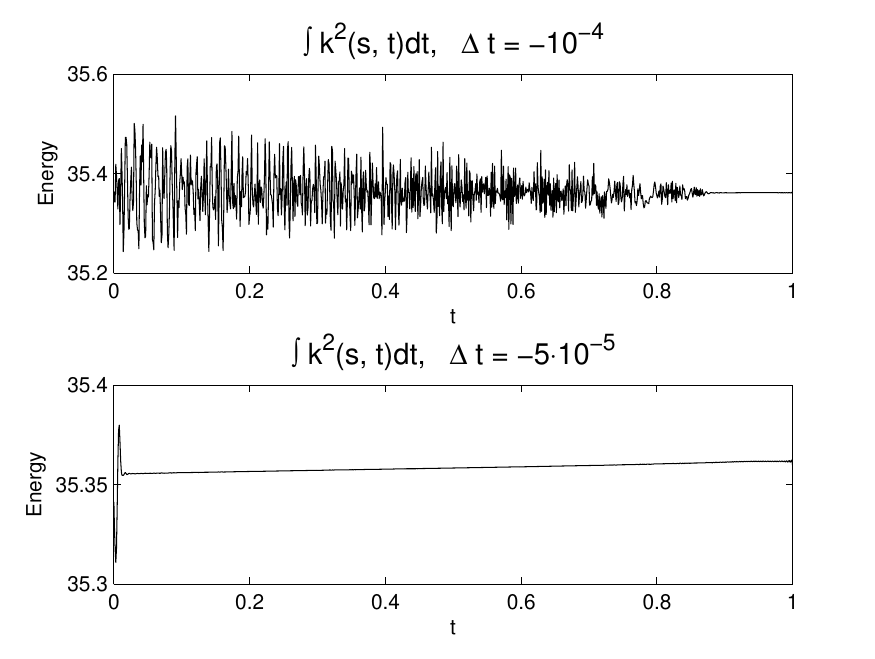}}
\caption{Numerical evolution of the energy in $s\in[s_a, s_b]$, with
initial data from \eqref{experimento1}. Unlike in the non-truncated
case with $s\in\mathbb R$, the energy in the interval $s\in[s_a,
s_b]$ keeps approximately constant .} \label{figura energy}
\end{figure}
at $s = 0$, the amount of energy \eqref{Esasb} in $s\in[s_a, s_b]$
grows up as we approach $t = 0$, tending to infinity. On the
contrary, in our numerical experiments, we have observed that the
energy in the interval $s\in[s_a, s_b]$ keeps approximately constant
after having fixed $\theta(s_a, t) = \theta^-$ and $\theta(s_b, t) =
\theta^+$; hence, we are preventing the energy outside the interval
$s\in[s_a, s_b]$ from entering.

In figure \ref{figura energy}, we have plotted the energy \eqref{Esasb}
in $s\in[s_a, s_b]$ for $\Delta t = -10^{-4}$ and $\Delta t = -5\cdot
10^{-5}$; remark that the energy conservation improves by diminishing
$|\Delta t|$.

On the other hand, the accuracy of $k(0, t)$ is rather poor for
small $t$ (see figure \ref{figura k(0,t)4096}).
\begin{figure}[!ht]
\centering \resizebox{9cm}{5cm}{\includegraphics{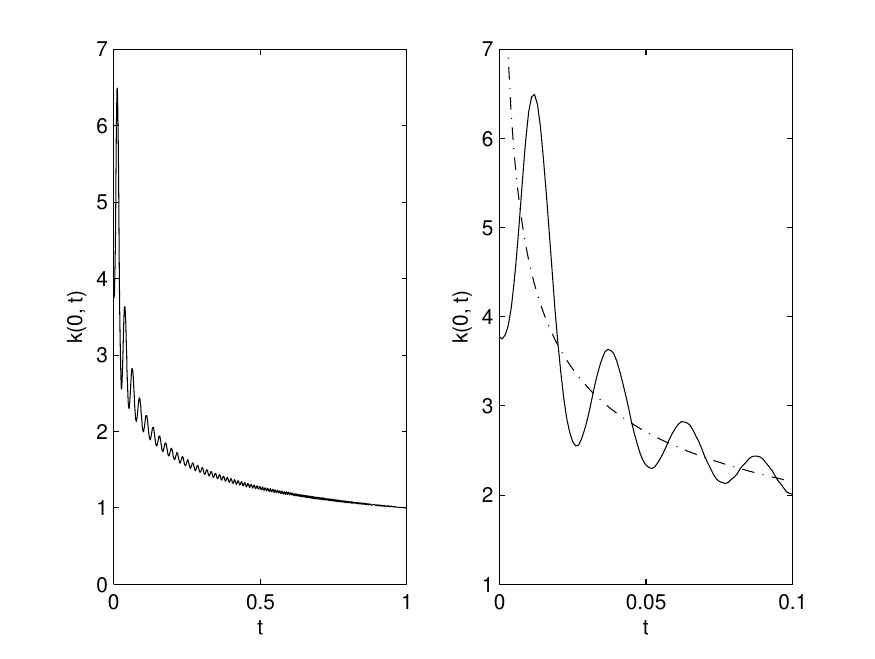}}
\caption{Numerical curvature at $s = 0$, with initial data from
\eqref{experimento1} and $\Delta t = -5\cdot10^{-5}$; the right-hand
side figure is a magnification of the left-hand figure; the exact
value \eqref{k0t} is plotted with a dotted stroke.} \label{figura
k(0,t)4096}
\end{figure}
Numerical experiments show that the way to recover the curvature $k(0,
t)$ with bigger accuracy near $t = 0$ is not to diminish $|\Delta t|$
or $\Delta s$, but rather to introduce more energy to the system, i.e.,
to lengthen the support of the initial datum $k(s, 1)$. To illustrate
it, we have executed the method also with the following initial data:
\begin{align}
\label{experimento2} & \begin{cases}
u(0) = 0.72, \qquad u'(0) = 1.1601860809647328, \\
u(x) \mbox{ integrated when $x$ in } [-400, 9.203], \\
s \in [-872.27, 309.22], \qquad N = 16384, \qquad \Delta s =
0.07211
\end{cases}
\intertext{and} \label{experimento3} & \begin{cases}
u(0) = 0.72, \qquad u'(0) = 1.1601860809647328, \\
u(x) \mbox{ integrated when $x$ in } [-800, 9.203], \\
s \in [-1744.62, 618.36], \qquad N = 32768, \qquad \Delta s =
0.07211.
\end{cases}
\end{align}

\begin{figure}[!ht]
\centering \resizebox{9cm}{5cm}{\includegraphics{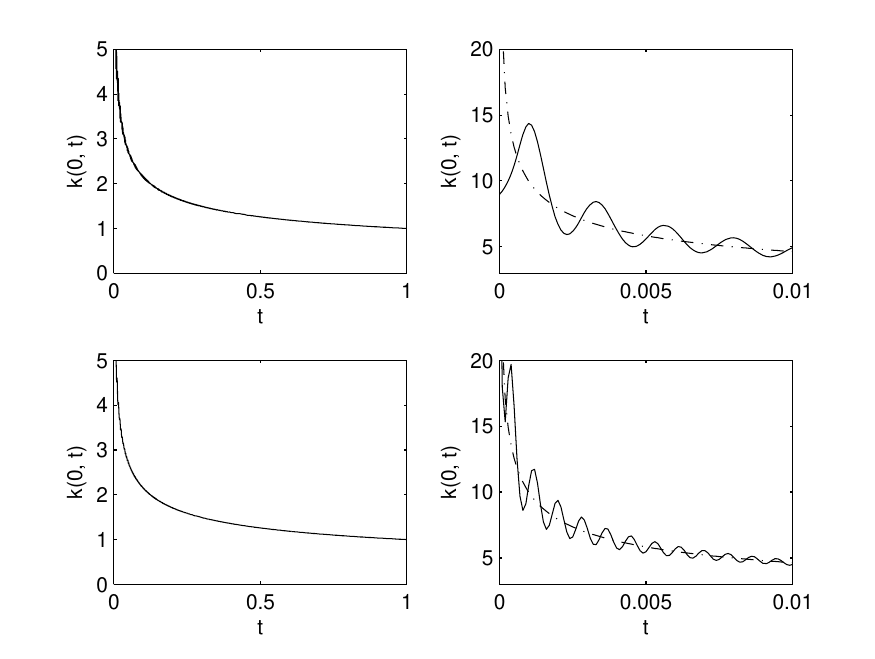}}
\parbox{14cm}{\caption{Curvature in $s = 0$. The upper part graphics
correspond to the experiment with data \eqref{experimento2}; the
lower part graphics correspond to the experiment with data
\eqref{experimento3}. The right-hand side graphics show $t\in[0,
0.01]$, i.e., time values close to the singularity; the exact value
\eqref{k0t} is plotted with dotted line.} \label{figura k(0,t)23}}
\end{figure}
\noindent Thus, we have made the initial support of $k(s, 1)$
approximately five and ten times as big respectively, but with the
same $\Delta s$ as in the previous experiment.

\begin{figure}[!ht]
\centering \resizebox{9cm}{5cm}{\includegraphics{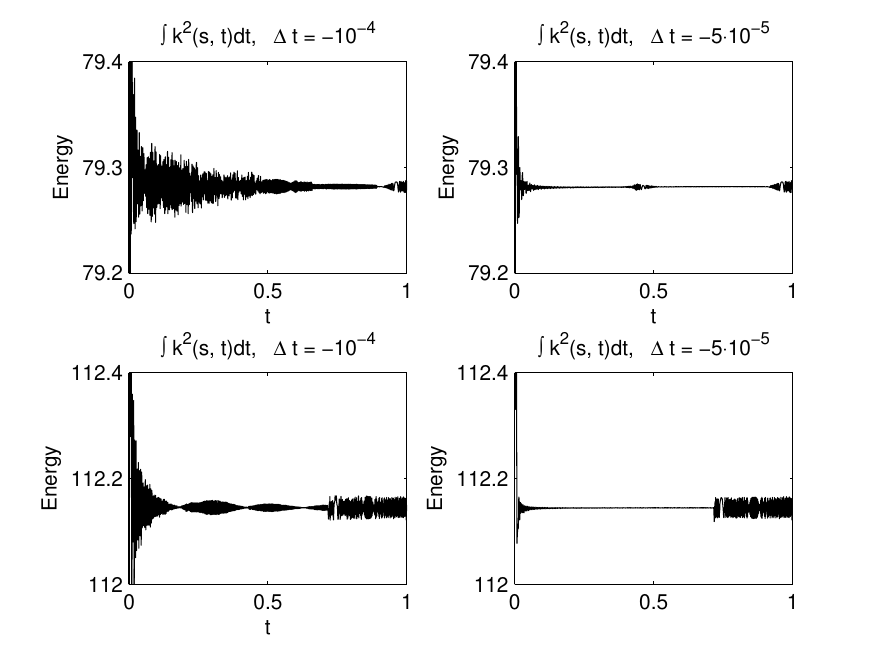}}
\caption{Numerical evolution of the energy. The upper part graphics
correspond to the experiment with data \eqref{experimento2}. the
lower part graphics correspond to the experiment with data
\eqref{experimento3}.} \label{figura energia16384}
\end{figure}

As we see in figure \ref{figura k(0,t)23}, we have been able to recover
$k(0, t)$ until much smaller $t$. Remark that the results are better in
the case where the support of $k(s, 1)$ is bigger (lower right-hand
side box). On the other hand, the observations about the energy are
valid again (see figure \ref{figura energia16384}).

\begin{figure}[!ht]
\centering \resizebox{9cm}{5cm}{\includegraphics{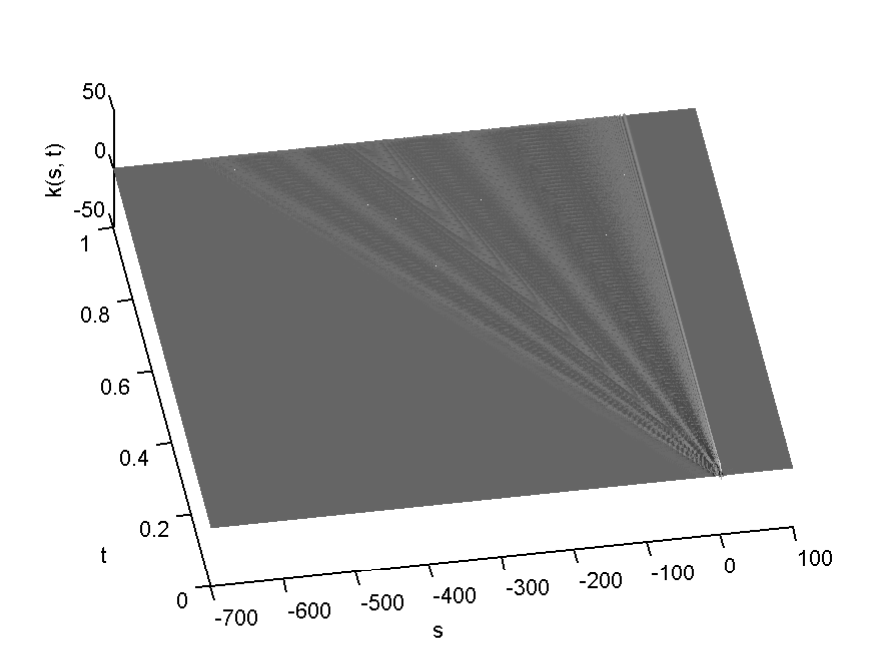}}
\caption{Numerical evolution of $k(s, t)$. The support of $k(s, t)$
tends to concentrate into $s = 0$ with lineal velocity, as $t\to0$.}
\label{figura evolucionk3D}
\end{figure}

In figure \ref{figura evolucionk3D}, we show the evolution of $k(s,
t)$ in function of space and time, for the simulation with initial
data \eqref{experimento2}. Clearly, the support of $k(s, t)$ tends
to concentrate into $s = 0$ with linear velocity, as $t\to0$.

Once we have obtained $\theta(s, t)$, it is immediate to recover $z(s,
t)$, from $\theta(s, t)$, making $z_s(s, t) = \exp(i\theta(s, t))$,
except for a rigid movement determined by \eqref{z(0,t)0}.
Nevertheless, from a numerical point of view, it is better for small
times to fix $z(s_b, t)\equiv z(s_b, 0)$. In figure \ref{figura3
z(s,t)}, we have superimposed the graphs of $z(s, t)$ in a neighborhood
of $s = 0$, with initial data \eqref{experimento2}, at times $t = 1$,
$t = 0.5$, $t = 0.1$, $t = 0.01$ and $t = 0.001$, having fixed $z_s(0,
t) = 1$. We can clearly appreciate the self-similar character of the
solutions.
\begin{figure}[!ht]
\centering \resizebox{9cm}{5cm}{\includegraphics{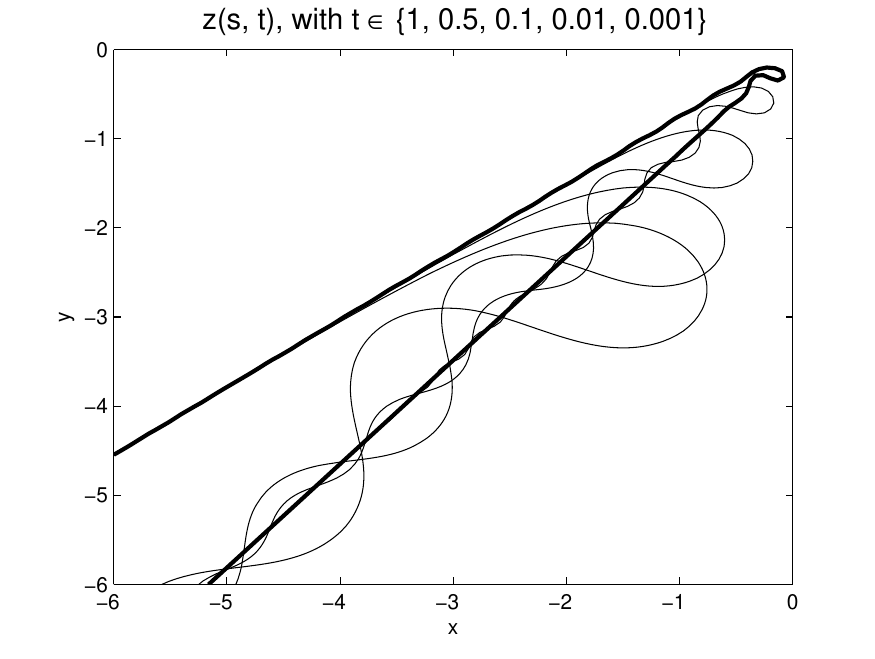}}
\caption{Numerical evolution of $z(s, t)$ for different times, with
initial data \eqref{experimento2}. The self-similarity is patent.
The thick dark curve corresponds to $t = 0.001$.} \label{figura3
z(s,t)}
\end{figure}

In figure \ref{figura3 pares sin autointersecciones}, we have plotted
with a thicker stroke the admissible pairs for which $z(s, t)$ has no
self-intersections, highlighting with a circle the pair that we have
used in our experiments. Figure \ref{figura3 z(s,t)} and figure
\ref{figura3 pares sin autointersecciones} explain why we have chosen
$(u(0), u'(0)) = (0.72, 1.1601860809647328)$. Indeed, we have taken an
admissible pair for which $z(s, t)$ has no self-intersections, but such
that it is near the situation when self-intersections happen; because
of that, $\int k(s, t)ds > 2.92$, i.e., a value near $\pi$. Theorems
\ref{teorema1} and \ref{teorema2} only guarantee the formation of the
singularity for small values of $\int_{-\infty}^{+\infty}k(s', t)ds' =
\theta^+ - \theta^-$, but the previous results and figure \ref{figura3
pares sin autointersecciones} give evidence that they also true for any
value of $\theta^+ - \theta^-\in(-\pi+\varepsilon, \pi-\varepsilon)$,
$0<\varepsilon\ll1$ (cf. section \ref{sectionintegral}). In any case,
all admissible data on the dark, 's' shaped curve in figure
\ref{figura3 pares sin autointersecciones} give self-similar corner
solutions.
\begin{figure}[!ht]
\centering
\resizebox{9cm}{5cm}{\includegraphics{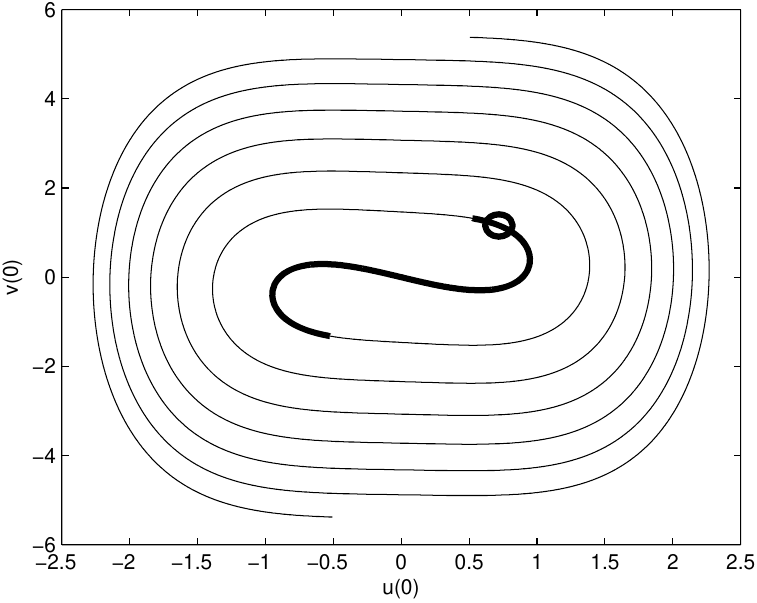}}
\caption{The thick dark curve shows the admissible pairs $(u(0),
u_x(0))$, such that $z(s, t)$ has no self-intersections. The circle
indicates the pair that we have considered.} \label{figura3 pares
sin autointersecciones}
\end{figure}

Before finishing this section, let us underline four main
conclusions:
\begin{itemize}

\item The finite energy in $[s_a, s_b]$ keeps approximately constant after
fixing the tangent vector of $z(s, t)$ at $s = s_a$ and $s = s_b$.

\item The support of $k(s, t)$, and hence, the energy,
tend to concentrate into $s = 0$ with lineal velocity, as $t\to0$.

\item $k(0, t)$ is well recovered, even for small $t$. To recover
$k(0, t)$ for smaller times, it is necessary to consider a bigger
support of $k(s, 1)$.

\item Our numerical results generalize theorems \ref{teorema1} and
    \ref{teorema1}, because they give evidence that the formation
    of the corner singularity happens for any value of $\theta^+ -
    \theta^-\in(-\pi+\varepsilon, \pi-\varepsilon)$,
    $0<\varepsilon\ll1$.

\end{itemize}

\noindent Hence, we are approximating a Dirac delta function
\eqref{ks0}, but we will never been able to create a singularity,
because the energy is finite. In the non-truncated problem, with
$s\in\mathbb R$, the energy is infinite for all $t$ and it tends to
concentrate at $s = 0$ as $t \to 0$, which causes the singularity
\eqref{ks0}.

Because of the parallelisms between the truncated and the non-truncated
problem, we can say that our numerical solutions reproduce the
non-truncated problem from a qualitative point of view, suggesting that
we could approach the Dirac delta function as much as we want, by
increasing the initial support of $k(s, 1)$. It would be very
interesting to prove analytically that we can recover the exact problem
by making $s_a\to-\infty$, $s_b\to+\infty$.

\section{$\int k(s', t)dt$, in terms of $u(0)$ and $u_x(0)$}

\label{sectionintegral}

As seen in section \ref{subsection ths1}, obtaining $\theta(s, 1)$ for
our experiments involves an estimate of $\int_{-\infty}^{+\infty}k(s',
1)ds'$. This integrate is hence determined by one admissible pair of
initial data $(u(0), u_x(0))$ of (eq.~\ref{airymodmu0}), i.e., one
point of the curve in figure \ref{figura u(x) pares admisibles}.

We have integrated \eqref{uxvxgx} in $x\in[-80, 20]$, $|\Delta x| =
10^{-5}$, for a large set of admissible pairs $(u(0), u_x(0))$ taken as
initial conditions. Following exactly the same steps as in
\ref{subsection ths1}, we have estimated $\int_{-\infty}^{+\infty}k(s',
1)ds'$.

\begin{figure}[!ht]
\centering \resizebox{9cm}{5cm}{\includegraphics{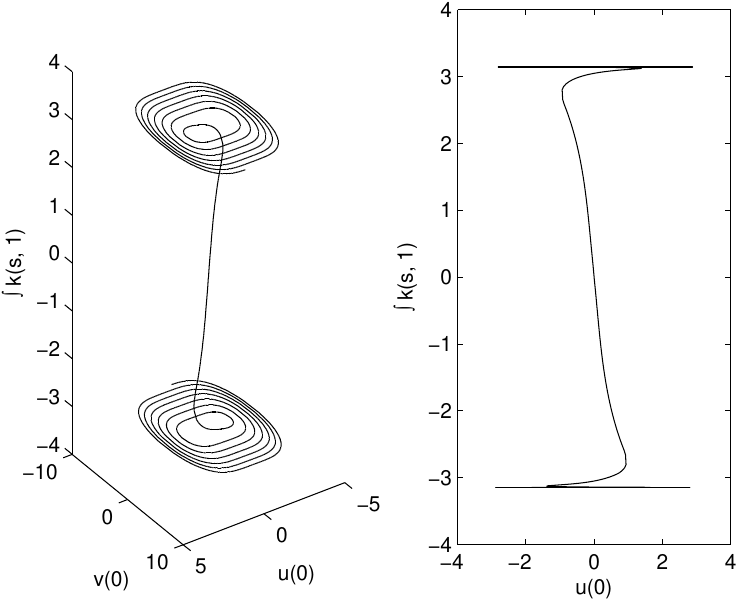}}
\caption{Integral of the curvature. The left-hand side shows the
admissible pairs $(u(0), u_x(0))$, together with their corresponding
integral $\int_{-\infty}^{+\infty}k(s', 1)ds'$. The right-hand side
is a side view of the left-hand side.} \label{figura3 integralk}
\end{figure}

In figure \ref{figura3 integralk}, each admissible pair $(u(0),
u_x(0))$ appears together with its corresponding estimate of
$\int_{-\infty}^{+\infty}k(s', 1)ds'$;
the curve we have plotted in this way is antisymmetric with respect to
the origin.

For the set of admissible data considered, our numerical estimates
satisfy
\begin{equation}
-\pi -5.6\cdot10^{-4} < \int_{-\infty}^{+\infty} k(s', 1, u(0), u'(0))ds'
< \pi+5.6\cdot10^{-4}.
\end{equation}

\noindent Hence, there is strong numerical evidence to accept that
\begin{equation}
\int_{-\infty}^{+\infty} u(x)dx \in\left(-\dfrac{\pi}{2},
\dfrac{\pi}{2}\right) \Longleftrightarrow \int_{-\infty}^{+\infty}
k(s', t)ds'\in(-\pi, \pi), \qquad\forall t.
\end{equation}

\noindent For values of the integral in $[-3-\varepsilon,
3+\varepsilon]$, $0 <\varepsilon \ll 1$, it is evident from the
right-hand side part, which represents the first and third
components of the curve on the left-hand side, that there is only
one corresponding initial condition for \eqref{airymodmu0}. The
values of the integral seem to converge exponentially to $\pm\pi$.
We are prone to think that there is one single initial condition for
each value in $(-\pi, \pi)$, although if we have exponential
convergence, it will be much more difficult to give numerical
evidence for the values nearest to $\pm\pi$.

\section{Using a simple closed $z(s, 1)$ as initial data}

\label{sectionclosed}

Coming back to our second experiment, with initial data
\eqref{experimento2}, and motivated by the vortex patch problem, it is
interesting to see what happens if we add a big loop to $z(s, 1)$, so
that it becomes a simple closed curve without intersections. For that
purpose, we have appended a smooth function to $\theta(s, 1)$ at $s =
s_a$, such that the new value for $\theta^+ - \theta^-$ is $\theta^+ -
\theta^-=2\pi$. This function is a rescaling of
\begin{align}
\Psi(s) & = \int_0^s\psi(s')ds', \quad s\in[0, 1], \\
\intertext{with} \psi(s) & =
\begin{cases}
0, & s\in [0, \alpha) \cup \{1\},
\\
\exp\left[\dfrac{-\beta}{(x-\alpha)(1 - x)}\right], & s\in(\alpha,
1),
\end{cases}
\end{align}

\noindent being both $\psi$ and $\Psi$ are regular in $[0, 1]$.

In what follows, we have rotated $z$ in such a way that $z_s(s_a) =
1$.
\begin{figure}[!ht]
\centering \resizebox{8cm}{4cm}{\includegraphics{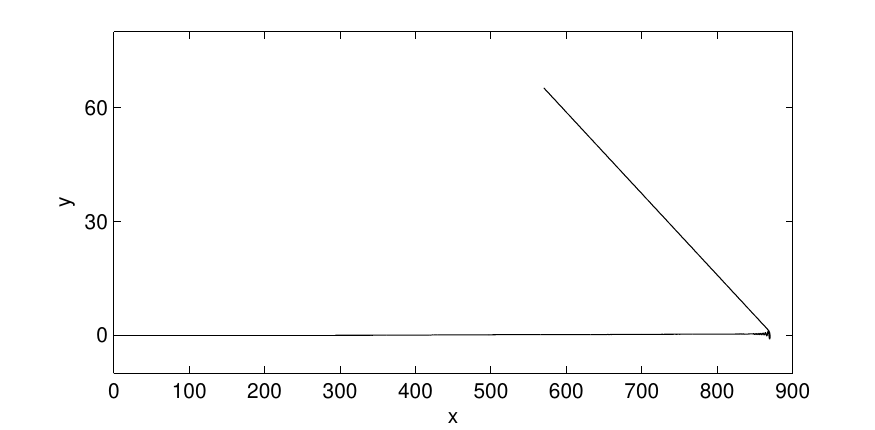}}
\caption{$z(s, 1)$ corresponding to the initial data
\eqref{experimento2} before being closed} \label{figura sinalargar}
\end{figure}
If we plot the $z$ corresponding to initial data \eqref{experimento2},
we observe in figure \ref{figura sinalargar} that the lower asymptotic
line of $z$ is much longer than the upper one. That is why we have
defined $\psi(s)$ being $\psi(s)\equiv 0$ in a certain $[0, \alpha]$,
in order to prolong the upper asymptotic line of $z$ before closing $z$
by means of a loop; the choice of $\alpha$ determines the additional
length that we add to the upper asymptotic line, whereas the choice of
$\beta$ determines the length of the loop. We have to adjust both
parameters wisely, so that $z$ becomes smoothly closed.

Remembering that $\theta(s, 1)$ was defined in such a way that
$\theta(s_b) = 0$, the new expression for $\theta(s, 1)$ looks like
this:
\begin{align}
\theta(s, 1) =
\begin{cases}
\theta(s, 1), & s\in[s_a, s_b]
\\
\dfrac{2\pi - \theta(s_a, 1)}{\Psi(1)}\Psi\left(\dfrac{s-s_b}{3(s_b
- s_a)}\right), & s\in[s_b, s_b + 3(s_b - s_a)],
\end{cases}
\end{align}

\noindent i.e., we have made four times as big the length of the new
$[s_a, s_b]$ and our new $s_b$ is $s_b = s_b + 3(s_b - s_a)$; the new
amount of nodes $s_j$ is $4N + 1 = 65537$, so that $\Delta s$ remains
unchanged. Notice that, with some notational abuse, $\theta(s, 1)$
refers to both the old and the prolonged functions and $s_b$ to the old
and the new right extreme. The new $\theta$ satisfies $\theta(s_b, 1) -
\theta(s_a, 1) = 2\pi$. This is a necessary but not sufficient
condition in order that $z$ may be closed. To determine $\alpha$ and
$\beta$, we proceed as follows:

If $\alpha$ is fixed, by a bisection technique, we determine a
$\beta$ such that $\Im(z(s_b, 1)) = 0$. Then, for some values of
$\alpha$, $\Re(z(s_b), 1)$ will be positive, while for some other
values, it will be negative. For instance, for two choices of
$\alpha$, $\alpha = 0.1$ and $\alpha = 0.2$, we have applied the
bisection technique to find the corresponding $\beta$ in $[0.5,
20]$. At the end of the process, $z(s_b, 1)$ is not exactly real,
but its imaginary part is negligible
\begin{align*}
\alpha & = 0.1 & \Rightarrow & & &
\begin{cases}
\beta = 5.492976875712412
\\
z(s_b, 1) \approx 209 + 5i\cdot10^{-12}
\end{cases}
\\
\alpha & = 0.2 & \Rightarrow & & &
\begin{cases}
\beta = 2.824543846424351
\\
z(s_b, 1) \approx -146 - i\cdot10^{-12}.
\end{cases}
\end{align*}

\noindent Now, we have to apply the bisection technique to $\alpha$.
For example, for $\alpha = 0.15$, we have $\beta =
3.970063345682464$ and $z(s_b, 1) \approx 32 - 2i\cdot10^{-12}$.
\begin{figure}[!ht]
\centering \resizebox{9cm}{5cm}{\includegraphics{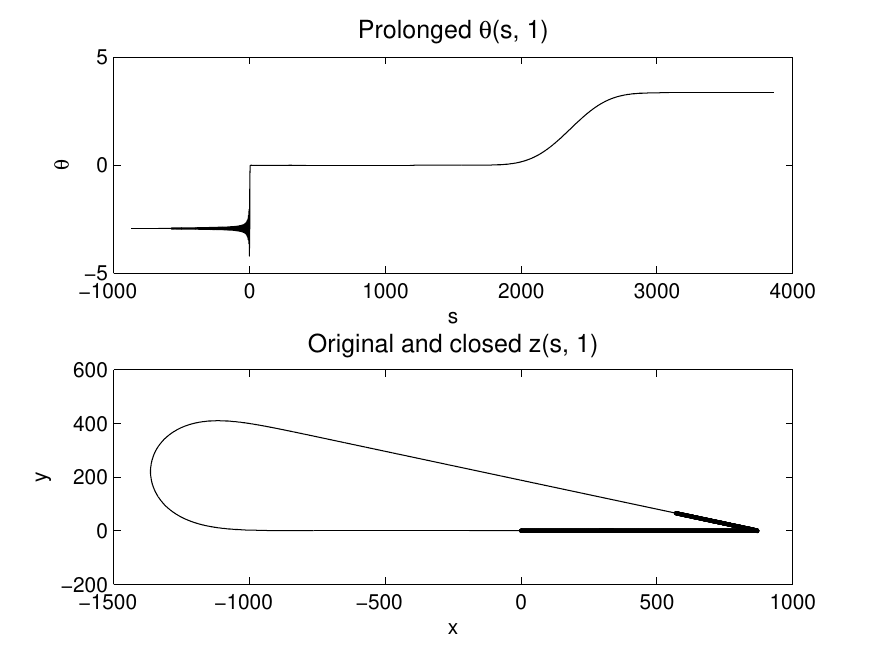}}
\caption{New $\theta(s, 1)$ and $z(s, 1)$} \label{figura zcerrada}
\end{figure}
That means that the correct $\alpha\in[0.15, 0.2]$. After several
iterations, the double-bisection algorithm yields
\begin{align}
\begin{cases}
(\alpha, \beta) = (0.158902181218767, 3.741222383167766)
\\
z(s_b, 1) \approx (2 + 3i)\cdot10^{-12}.
\end{cases}
\end{align}

\noindent Since $s\in[-872.27, 3853.69]$ and $z(s, 1)$ is
parameterized respect to the arc-length parameter, the length of the
curve is $4725.96$.  In figure \ref{figura zcerrada}, we can see the
prolonged version of $\theta(s, 1)$, as well as $z(s, 1)$, with and
without the loop. We have closed $z(s, 1)$ almost in a perfect way,
because $z(s_b, 1) - z(s_a, 1) = \mathcal O(10^{-12})$.

When $z$ is closed, $z_s(0, t)$ is obviously no longer constant and
\eqref{z(0,t)0} does not hold any more; in fact, $z(s, t)$ will
rotate, so we need to give the evolution of a point and its angle
for all $t$, which is quite straightforward. We just have to
integrate the following two ODE's in $t$ for some $s_0$:
\begin{align*}
\theta_t(s_0, t) & = -\theta_{sss}(s_0, t) -
\dfrac{1}{2}\theta_s^3(s_0, t),
\\
z_t(s_0, t) & = - z_{sss}(s_0, t) + \dfrac{3}{2}\bar z_s(s_0, t)
z_{ss}^2(s_0, t).
\end{align*}

\noindent The right-hand side of the first equation is known, so we
can update $\theta(s_0, t^n)$. Then, bearing in mind that $z_s(s_0,
t^n) = \exp(i\theta(s_0, t^n))$, the right-hand side of the second
equation is also known and we update $z(s_0, t^n)$.

In any case, since this paper aims at illustrating the energy
concentration process, we will not rotate $z$, in order to compare
the situations with and without the loop.

\subsection{Numerical experiments}

\label{subsection experimentsthclosed}

\begin{figure}[!ht]
\centering \resizebox{9cm}{4cm}{\includegraphics{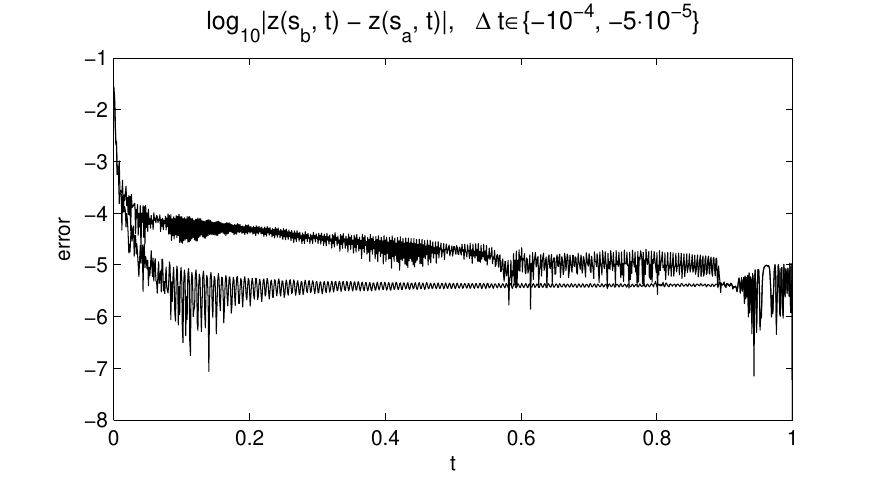}}
\caption{Error when closing $z(s, t)$. The error is smaller for $\Delta t = -5\cdot 10^{-5}$.} \label{figura errorz}
\end{figure}

We have executed the method with the new $\theta(s, 1)$. Remark that
the evolution variable is $\theta$, not $z$, so the method does not
guarantee a priori that $z(s_a, t) = z(s_b, t)$, for all $t$.
\begin{figure}[!ht]
\centering \resizebox{9cm}{5cm}{\includegraphics{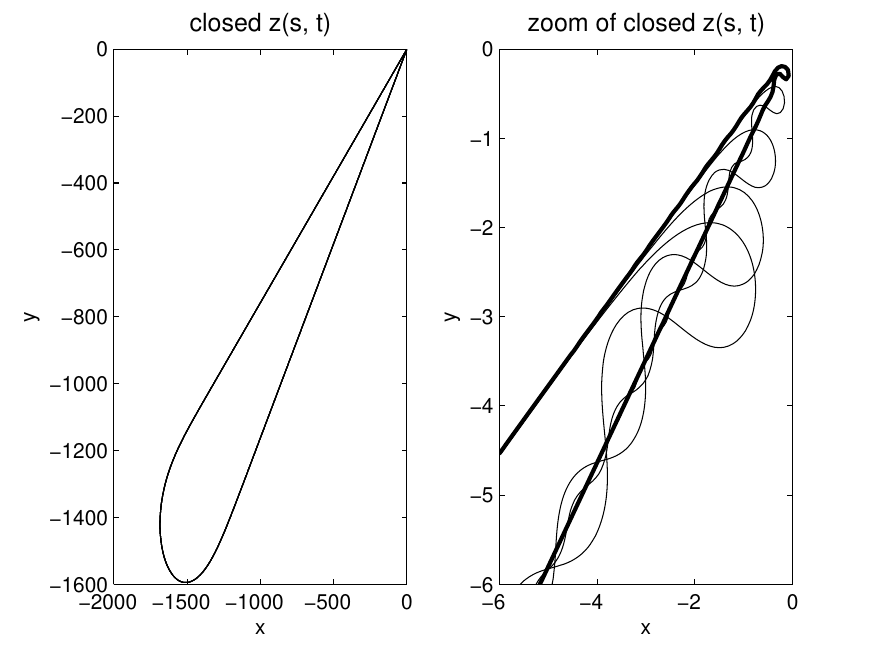}}
\caption{$z(s, t)$ and zoom of $z(s, t)$, for $t\in\{1, 0.5, 0.1,
0.01, 0.001\}$. The thick dark curve at the right-hand side
corresponds to $t = 0.001$.} \label{zevol}
\end{figure}
Nevertheless, the curve keeps closed with great accuracy, as we can
see in figure \ref{figura errorz}, where we have plotted the decimal
logarithm of $|z(s_b, t) - z(s_a, t)|$, for $\Delta t$ small enough.
Indeed, for $\Delta t = -5\cdot10^{-5}$, $|z(s_a, t) - z(s_b, t)| <
10^{-4}$ until very small times and, since the length of the curve
is approximately $4726$, that means a relative error smaller than
$10^{-7}$.

On the other hand, the loop does not have any visible influence on
the formation of the corner-shaped singularity for $z(s, t)$. In
figure \ref{zevol}, we have plotted $z(s, t)$, with $z_s(0, t) = 1$,
for $t = 1$, $0.5$, $0.1$, $0.01$ and $0.001$. In the left-hand side
part, all the plots seem to be a single one; nevertheless, after a
potent zoom, the formation of the singularity at $s = 0$ is clearly
exhibited. In fact, the right-hand side of the figure is not
ocularly distinguishable from figure \ref{figura3 z(s,t)}.

From the left-hand side of figure \ref{zevol}, it is obvious that
the enclosed area remains almost constant, for all $t$. We can check
this by means of the well known formula:
\begin{align}
\mbox{Area} = \dfrac{1}{2}\oint(xdy - ydy) =
\dfrac{1}{2}\int_{s_a}^{s_b}[x(s)y'(s) - x'(s)y(s)]ds, \quad z = x +
iy,
\end{align}

\noindent which is a direct application of Green's theorem. Figure
\ref{figura zarea} shows that the enclosed area is preserved with
great accuracy.
\begin{figure}[!ht]
\centering \resizebox{9cm}{4cm}{\includegraphics{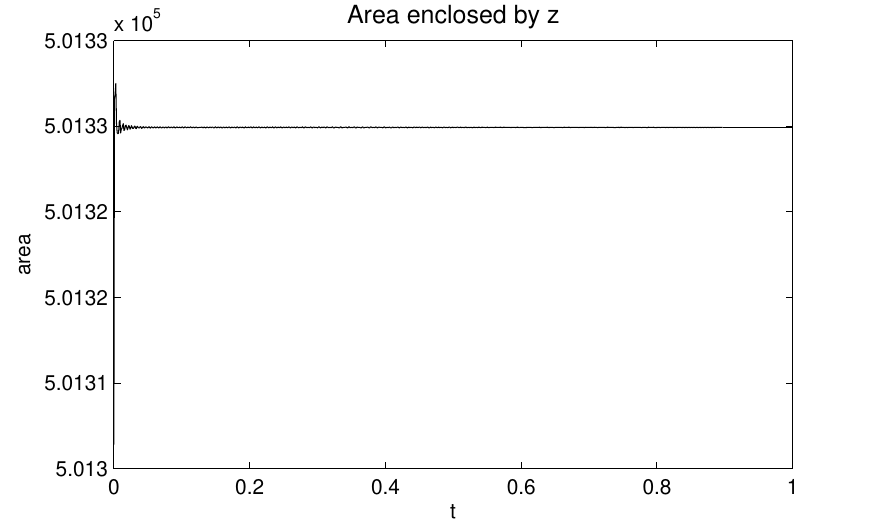}}
\caption{Preservation of the area enclosed by $z$} \label{figura
zarea}
\end{figure}

\section{Conclusions}

\label{sectionconclusions}

In this paper, we have consider a geometric planar flow
(eq.~\ref{flujo}) found by Goldstein and Petrich \cite{goldstein},
while considering the evolution of a vortex patch. Perelman and Vega
proved in \cite{perelman} that it has a one-parameter family of
solutions that develop corner-shaped singularities at finite time. We
have studied those solutions from a numerical point of view, trying to
reproduce the singularity formation.

The flow can be characterized by the equation of its curvature, $k$,
which is the modified Korteweg-de~Vries \eqref{mkdv} equation.
Nevertheless, we have found the angle $\theta = \int^s kds'$ to be the
most adequate evolution variable, because it allows to preserve
naturally the arc-length parameter and the angle formed by the two
asymptotic lines of $z$. Due to the difficulty of considering the whole
$\mathbb R$, we have taken $s\in[s_a, s_b]$. The main result is that,
even after fixing $\theta(s_a, t) = \theta^-$, $\theta(s_b, t) =
\theta^+$, i.e., after fixing the tangent vectors of $z$ at $s_a$ and
$s_b$, we are still able to approximate the formation of the
corner-shaped singularity. Indeed, the energy
$\int_{s_a}^{s_b}k^2(s',t)ds' < \infty$ keeps approximately constant
and it tends to concentrate at $s = 0$; moreover, we recover with great
accuracy $k(0, t)$, even for small $t$, tending $k$ to a Dirac delta
function. The numerical results suggest that the accuracy of $k(0, t)$
could be improved arbitrarily by increasing the length of $[s_a, s_b]$;
it would be very interesting to prove analytically that we can recover
the solution of the exact problem by making $s_a\to-\infty$ and
$s_b\to\infty$.

The process of obtaining the initial $\theta$ involves an estimate
of the conserved quantity $\int_{-\infty}^{+\infty} k(s', t)ds'$. We
have given numerical evidence that $|\int_{-\infty}^{+\infty} k(s',
t)ds'| < \pi$.

Motivated by the vortex patch problem, we have also considered a
regular simple closed curve as initial datum, by appending a regular
function to $\theta(s, 1)$ in such a way that $z(s, 1)$ is closed by
a big loop. The method preserves the area enclosed by the curve with
great accuracy and, more interestingly, the loop appears to have no
influence on the energy concentration process; nevertheless, no
singularities can be expected because the energy of our curves is
finite \cite{kenig,bourgain}. This is in agreement with the vortex
patch theory, where it is known that no singularities may arise from
smooth initial contours \cite{chemin, bertozzi}. In \cite{perelman},
Perelman and Vega considered instead infinite energy solutions of
the mKdV equation, which causes the singularity to happen.

One of the shortcomings of our model is that the arc-length
parameter is preserved for a simple closed curve and, hence, the
total length of the curve remains constant. In the vortex patch
problem, on the contrary, there are examples in which the length or
curvature of the vortex patch boundary grow rapidly \cite{alinhac,
constantin}. Therefore, the mKdV equation should be complemented by
another equation for the evolution of $s_\alpha$.

\section{Acknowledgements}

The author would want to express his gratitude to L.~Vega and C.
García-Cervera for very valuable advice concerning this paper.

\bibliographystyle{siam}

\end{document}